\RequirePackage{fix-cm}
\documentclass[smallextended]{svjour3}       
\smartqed  

\usepackage{etoolbox}
\newcommand*{\affaddr}[1]{#1} 
\newcommand*{\affmark}[1][*]{\textsuperscript{#1}}

\usepackage{amssymb}
\usepackage{amsmath,bm}
\usepackage{mathabx}
\usepackage{graphicx}
\usepackage{eurosym}
\usepackage{enumitem}
\usepackage{hyperref}
\usepackage{booktabs}
\usepackage{multirow}
\usepackage{siunitx}
\usepackage{textgreek}
\usepackage{color}
\usepackage[lined,ruled]{algorithm2e}
\usepackage[caption=false]{subfig}
\newcommand{\plus}{\scalebox{.60}{$+$}}
\newcommand{\minus}{\scalebox{.60}{$-$}}

\newtheorem{cond}{Condition}
\newtheorem{proper}{Property}

\makeatletter
\newcommand*{\algotitle}[2]{%
  \stepcounter{algocf}%
  \hypertarget{algocf.title.\theHalgocf}{}%
  \NR@gettitle{#1}%
  \label{#2}%
  \addtocounter{algocf}{-1}%
}
\makeatother

\begin{document}


\title{Adaptive Benders decomposition and enhanced SDDP for multistage stochastic programs with block-separable multistage recourse
\thanks{Research is supported by the Engineering and Physical Sciences Research Council (EPSRC) through the CESI project (EP/P001173/1)}}


\author{Nicol\`o Mazzi\protect \and Ken McKinnon \and Hongyu Zhang\thanks{Corresponding author \email{hongyu.zhang@soton.ac.uk}}}


\institute{N. Mazzi\affmark[1] \at
           \email{nicolo.mazzi@ed.ac.uk} \and
           K. McKinnon\affmark[1] \at
           \email{K.McKinnon@ed.ac.uk} \and 
           H. Zhang\affmark[2] \at
           \email{hongyu.zhang@soton.ac.uk} \and 
           \affaddr{\affmark[1] School of Mathematics, University of Edinburgh, James Clerk Maxwell Building, Edinburgh, EH9 3FD, United Kingdom}\\
            \affaddr{\affmark[2] School of Mathematical Sciences, University of Southampton, Building 54, Highfield Campus, Southampton, SO14 3ZH, United Kingdom}}

\date{version of \textbf{\today}} 

\maketitle

\begin{abstract}
This paper proposes an algorithm to efficiently solve multistage stochastic programs with block separable recourse where each recourse problem is a multistage stochastic program with stage-wise independent uncertainty. The algorithm first decomposes the full problem into a reduced master problem and subproblems using Adaptive Benders decomposition. The subproblems are then solved by an enhanced SDDP. The enhancement includes (1) valid bounds at each iteration, (2) a path exploration rule, (3) cut sharing among subproblems, and (4) guaranteed $\delta$-optimal convergence. The cuts for the subproblems are then shared by calling adaptive oracles. The key contribution of the paper is the first algorithm for solving this class of problems. The algorithm is demonstrated on a power system investment planning problem with multi-timescale uncertainty. The case study results show that (1) the proposed algorithm can efficiently solve this type of problem, (2) deterministic wind modelling underestimate the objective function, and (3) stochastic modelling of wind leads to different investment decisions.

\end{abstract}

\section{Introduction}\label{sec:Intro}
\begin{figure}[!htb]
    \centering
    \includegraphics[scale=0.8]{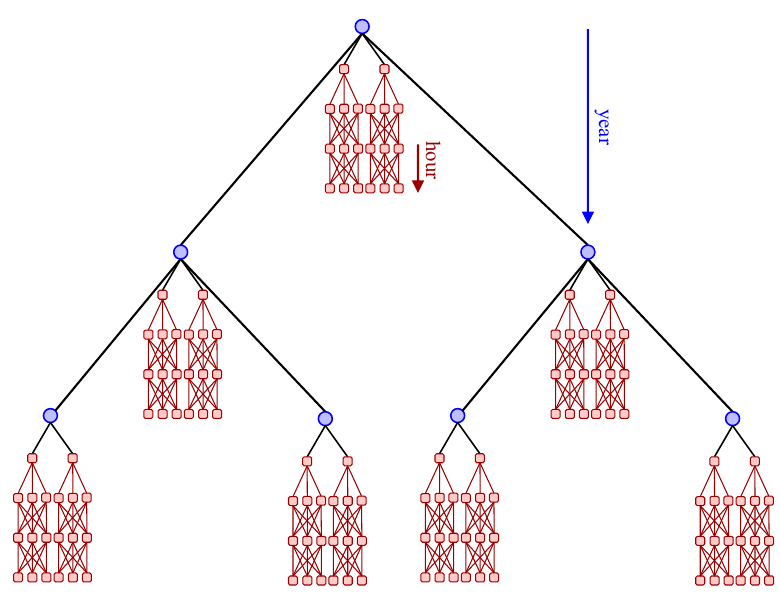}
    \caption{Illustration of a stochastic programming scenario tree with long-term uncertainty (blue circles) and short-term uncertainty (red squares).}
    \label{fig:long-term tree}
\end{figure}

\begin{figure}[htbp!] 
  \centering
      \includegraphics[width=1.00\columnwidth]{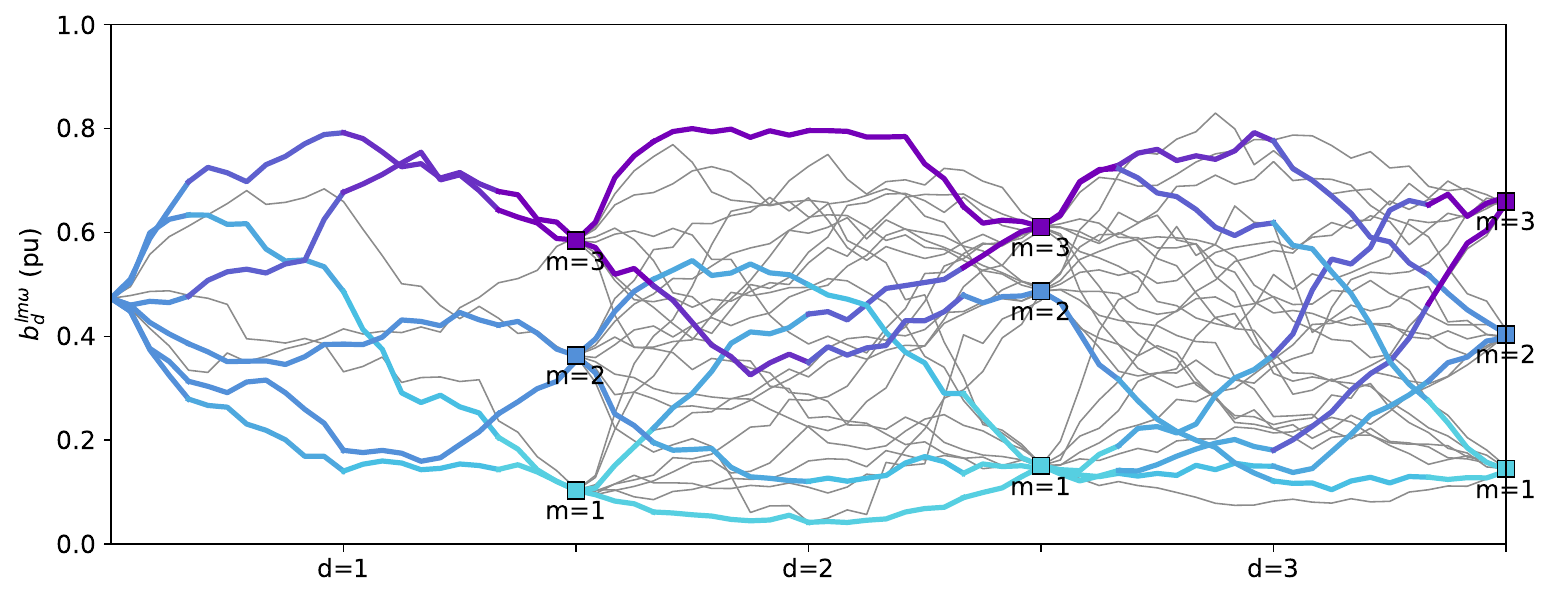}
	  \caption{Illustration of a short-term scenario folded tree (short-term wind power uncertainty for 3 days). The grey continuous lines represent the wind trajectories $b^{lm\omega}_d$, and the squares at the end of each day indicate the state $m$ they represent. Some of the trajectories are highlighted to show some of the possible paths.}
      \label{fig:short-term tree}
\end{figure}

In this paper, we are interested in solving problems of form
\begin{equation}\label{eq:opt_init}
	\begin{aligned}
		\mathbf{MP} : \hspace{15pt} \underset{\mathbf{x} \in \mathcal{X}}{\text{min}} & \;\; f(\mathbf{x}) + \sum_{i \in \mathcal{I}} \pi_i \hspace{1pt} g(x_i), 
	\end{aligned}
\end{equation}
where $\mathcal{X}$ is a compact set, $x_i$ is a (possibly overlapping) subvector of $\mathbf{x}$, and $\pi_i$ are non-negative coefficients. The value of $g(x)$ is obtained solving a multistage stochastic linear program with $\mathbf{d}$ stages, defined as
\begin{multline}\label{eq:multistage_g}
    \mathbf{SP}(x): \quad g(x) :=  \mathbb{E}_1 \hspace{-2pt} \left[ \hspace{-2pt} \underset{ \substack{ y_1 \in \mathcal{Y}_1 \\ A_1 y_{1}  \leq  x^{\top} \hspace{-2.pt} B_1 b_1}}{\text{min}} \hspace{-12pt} c^\top_1 y_1 \hspace{-1pt} + \hspace{-1pt} \mathbb{E}_{2|\tilde{b}_1} \hspace{-2pt} \left[ \hspace{-2pt} \underset{ \substack{ y_2 \in \mathcal{Y}_2 \\ A_2 y_2  \leq  x^{\top} \hspace{-2.pt} B_2 b_2 + C_2 y_1 }}{\text{min}} \hspace{-22pt} c^\top_2 y_2 + \; ... \; + \right. \right. \\
    \left. \left. + \mathbb{E}_{d|\tilde{b}_{d\minus1}} \hspace{-2pt} \left[ \hspace{-2pt} \underset{ \substack{ y_d \in \mathcal{Y}_d \\A_d y_d  \leq  x^{\top} \hspace{-2.pt} B_d b_d + C_d y_{d\minus1} }}{\text{min}} \hspace{-26pt} c^\top_d y_d  \hspace{-1pt} +  \; ... \; + \hspace{-1pt} \mathbb{E}_{\mathbf{d}|\tilde{b}_{\mathbf{d}\minus1}} \hspace{-2pt} \left[ \hspace{-2pt} \underset{ \substack{ y_\mathbf{d} \in \mathcal{Y}_\mathbf{d} \\A_\mathbf{d} y_\mathbf{d}  \leq  x^{\top} \hspace{-2.pt} B_\mathbf{d} b_\mathbf{d} + C_\mathbf{d} y_{\mathbf{d}\minus1} }}{\text{min}} \hspace{-26pt} c^\top_\mathbf{d} y_\mathbf{d} \hspace{4pt} \right] \hspace{-2pt}\right] \hspace{-2pt}\right] \hspace{-2pt}\right],
\end{multline}
where $A_d$, $B_d$, and $C_d$ are matrices of coefficients, $b_d$ is a random vector of coefficients, and $\mathcal{Y}_d$ is a set of linear constraints. Each stage $d$ is associated with a discrete set $\mathcal{M}$ of possible states and a discrete set $\Omega$ of possible scenarios, and $b^{lm\omega}_d$ is the realisation of $b_d$ that lands at state $m$ via scenario $\omega$ given that day $d-1$ ended at state $l$. Each realisation $b^{lm\omega}_d$ is associated to a probability $\pi^{lm}_d \times \pi^{\omega}$ of occurring, where $\pi^{lm}_d$ is the probability of landing in state $m$ given that day $d-1$ ended at state $l$ ($\sum_{m\in \mathcal{M}} \pi^{lm}_d = 1$) and $\pi^{\omega}$ is the probability of scenario $\omega$ which we assume equal for each scenario ($ \pi^{\omega} = \frac{1}{|\Omega|}, \,  \omega \in \Omega$). An example of such problem is illustrated in Figure \ref{fig:long-term tree}. An example of the short-term SDDP problem is illustrated in Figure \ref{fig:short-term tree}.

\subsection{Prior work}
\label{sec:prior work}

There has been no algorithm proposed to tackle this kind of problem. The main idea of the proposed algorithm is to first decompose the problem using Benders-type decomposition and then apply a modified version of Stochastic Dual Dynamic Programming (SDDP) to solve each block of multistage subproblem.

Multistage stochastic programming with block-separable recourse was proposed by \cite{Louveaux1986MultistageRecourse}. It was later reinvented by \cite{Kaut2014} from the standpoint of long-term infrastructure planning problems involving uncertainty from short-term and long-term time horizons, and the authors called it Multi-Horizon Stochastic Programming (MHSP). MHSP has then been applied to multiple energy system planning problems with short-term uncertainty \cite{Backe2022,Durakovic2024DecarbonizingSectors} and both short-term and long-term uncertainty \cite{Zhang2024AUncertainty}. However, there has been no study using MHSP and including short-term uncertainty that reveals multiple stages. In all existing literature using MHSP with short-term uncertainty, the short-term uncertainty is represented by some time slices of the short-term time series. Because of this, there are no recourse actions modelled, hence no flexibility. This is mainly due to computational limits. However, capturing flexibility in decision-making by having the option of delaying decisions is the key advantage and value of stochastic programming \cite{King2024ModelingEdition}. Therefore, one can argue that the way short-term uncertainty was modelled is not stochastic programming. In this paper, we aim to address this research gap by first modelling multistage short-term uncertainty and proposing Adaptive Benders decomposition and enhanced SDDP to solve the problem.

Benders decomposition was proposed to solve problems with complicating variables \cite{Benders1962}. It was applied to solve two-stage stochastic programming and was referred to as the L-Shaped method \cite{VanSlyke1969}. \cite{Birge1985DecompositionPrograms} proposed a nested Benders decomposition to solve multistage stochastic programmes. There have been several enhancements to Benders decomposition, such as stabilisation \cite{Zhang2024AUncertainty,Zhang2025IntegratedDecomposition}, cut selection, and inexact oracles \cite{Zakeri,Mazzi2020}. Benders decomposition can be applied to solve multistage stochastic programmes with block-separable recourse \cite{Zhang2024DecompositionProgramming}. Also, the block-separable structure can be exploited for efficient solution algorithms via cut sharing. The Adaptive Benders decomposition was proposed in \cite{Mazzi2020} to solve large-scale optimisation problems and is applicable for MHSP problems \cite{Zhang2024AUncertainty}. However, Adaptive Benders decomposition can suffer from oscillation like other types of Benders decomposition. Therefore, \cite{Zhang2024AUncertainty} proposed a level method stabilisation and achieved significant improvement, and the stabilised Adaptive Benders decomposition is able to solve linear programming problems with up to 1 billion variables and 4.5 billion constraints. A centred point stabilisation was proposed to solve problems with integer variables in the reduced master problem \cite{Zhang2025IntegratedDecomposition} and was applied to solve an integrated European energy system planning problem. Despite the development of Benders decomposition for MHSP, there is a research gap remaining unaddressed: there is no flexibility that can be captured in the block-separable subproblems. However, making the subproblem stochastic leads to a significant increase in computational difficulty. Therefore, in this paper, we propose to use Adaptive Benders decomposition to decompose the monolithic problem into a master problem and some blocks of subproblems. Then we use an enhanced SDDP to solve the subproblems.

SDDP was proposed by \cite{pereira1991} to solve hydropower scheduling problems under uncertainty. It is based on the assumption of uncertain parameters being stage-wise independent. In our paper, we also rely on this assumption for the short-term uncertainty. SDDP has been applied to many applications mainly in the energy field, such as \cite{Lara2020,Papavasiliou2018ApplicationUncertainty}. SDDiP was proposed to solve problems with integer variables in the subproblems \cite{Zou2019StochasticProgramming} and was applied to power system planning problems \cite{Lara2020}. \cite{Downward2020StochasticUncertainty} proposed an extended SDDP algorithm to solve problems with state-wise-dependent objective uncertainty. We refer to \cite{Fullner2021StochasticVariants} for a review of the SDDP algorithms and its variants. As a sampling-based algorithm, the original SDDP only has a statistical upper bound, assuming a minimisation problem. This makes the gap at each iteration invalid. The stopping criteria of the SDDP algorithm have also been criticised due to this \cite{Shapiro2011AnalysisMethod}. Also, there has been effort to improve the stopping criteria \cite{Homem-De-Mello2011SamplingScheduling,Leclere2020ExactDuality}. In this paper, we address this issue by proposing an enhanced SDDP which can obtain valid bounds at each iteration. In addition, we also propose a path search rule and cut sharing among subproblems. Finally, we can obtain guaranteed $\delta$-optimal convergence.

\subsection{Contributions}
\label{sec:contributions}

Our contribution is (1) Adaptive Benders decomposition with enhanced SDDP to solve block-separable multistage stochastic programmes with multistage recourse subproblems, (2) valid bounds for SDDP at each iteration, (3) a path exploration rule, (4) cut sharing among subproblems and convergence guarantee and proof of convergence of the enhanced SDDP.

\subsection{Paper structure}
\label{sec:paper structure}

The remainder of the paper is organised as follows. Section \ref{sec:ProbMod_Assump} introduces the assumptions and decomposition strategy. Section \ref{sec:sddp_approx} presents the SDDP algorithm with valid bounds and a path search rule. Section \ref{sec:sddp_extended} introduces the proposed enhanced SDDP algorithm with valid bounds, a path search rule and cut sharing. Section \ref{sec:Benders} presents the Adaptive Benders decomposition. Section \ref{sec:casestudy} tests the proposed algorithm on a power system investment planning problem. Finally, conclusions are drawn in Section \ref{sec:conclusions}.

\section{Problem modification and assumptions}\label{sec:ProbMod_Assump}
In the following, we introduce problem reformulation and assumptions.

\subsection{Problem modification}\label{subsec:ProbMod}
We separate the investment problem from the operational problems with a Benders decomposition approach. At iteration $w$, we solve the relaxation
\begin{equation}\label{eq:RMP}
	\begin{aligned}
		\mathbf{RMP} : \hspace{15pt} \underset{\mathbf{x} \in \mathcal{X}}{\text{min}} & \;\; f(\mathbf{x}) + \sum_{i \in \mathcal{I}} \pi_i \hspace{1pt} \beta_i, \\
		\text{s.t.} & \;\; \beta_i \geq \underline{\theta}^{\mathrm{w}}_{i} + \underline{\lambda}^{\mathrm{w} \top}_i \hspace{-1.5pt} \big( x_i - x^{\mathrm{w}}_i \big), \quad   \mathrm{w}=1,..,w-1, \;  i \in \mathcal{I},\\
	\end{aligned}
\end{equation}
of the $\mathbf{MP}$ problem in~\eqref{eq:opt_init}, which yields an optimal solution $\mathbf{x}^{w}$. To generate a new valid cutting plane $(\underline{\theta}^w_i,\underline{\lambda}^w_i)$ we solve subproblem $\mathbf{SP}$ in~\eqref{eq:multistage_g} and obtain a valid lower bound $\underline{\theta}^w_i$ on the true objective $\theta^w_i = g(x^w_i)$, and a valid subgradient $\underline{\lambda}^w_i$ w.r.t. $x_i$.
The master problem decisions $x_i$ are therefore fixed parameters for subproblem $\mathbf{SP}(x_i)$, and we reformulate~\eqref{eq:multistage_g} in a recursive fashion as $\theta^k_i:=\sum_{m \in \mathcal{M}} \pi^m_1 V^{m}_1$, where $V^{m}_1$ is the optimal solution of\\

\noindent $\mathrm{SP}^m_{1}$: \hspace{8pt} $V^m_1 :=$
\begin{equation}\label{eq:sp_init_0}
	\begin{aligned}
	\underset{y_1 \in \mathcal{Y}_1}{\text{min}} & \hspace{3pt} \sum_{\omega\in\Omega} \pi^{\omega} \left(  c^{\top}_1 y^\omega_1 + \sum_{n\in \mathcal{M}} \pi^{mn}_2 V^{mn}_{2} (y^\omega_1)\right)\\
    \text{s.t.} & \hspace{3pt} A_1 \hspace{1.0pt} y^\omega_{1} \leq x^{\top} \hspace{-1pt} B_1 \hspace{1.0pt} b^{m\omega}_{1}, \;  \omega \in \Omega, \\
	\end{aligned}
\end{equation}	
where $y_d = \{y^\omega_d, \omega \in \Omega \} $. Function ${V}^{lm}_{d}(y_{d\minus1})$ is defined as\\

\noindent $\mathrm{SP}^{lm}_{d}(y_{d\minus1})$: \hspace{8pt} $V^{lm}_d (y_{d\minus1}) :=$
\begin{equation}\label{eq:sp_init_d} 
	\begin{aligned}
	\underset{y_d \in \mathcal{Y}_d}{\text{min}} & \hspace{3pt} \sum_{\omega \in \Omega} \pi^{\omega} \hspace{-2pt} \left(c^{\top}_d y^\omega_d  + \hspace{-2pt} \sum_{n\in \mathcal{M}} \pi^{mn}_{d\plus1} V^{mn}_{d\plus1} (y^\omega_d) \right)\\
    \text{s.t.} & \hspace{3pt} A_d  \hspace{1pt} y^\omega_{d}  \leq  x^{\top} \hspace{-1.5pt} B_d \hspace{.5pt} b^{l m \omega}_d  + C_d \hspace{.5pt} y_{d\minus1}, \;  \omega \in \Omega. \\
	\end{aligned}
\end{equation}
Note that in the general case problem~\eqref{eq:sp_init_d} is not separable w.r.t. scenarios $\omega \in \Omega$ since $y_d \in \mathcal{Y}_d$ might include constraints (e.g., non anticipativity constraints) that links the scenarios together. At the last stage, the value ${V}^{lm}_{\mathbf{d}}(y_{\mathbf{d}\minus1})$ is obtained by imposing the cost-to-go function to 0, i.e.,\\ 

\noindent $\mathrm{SP}^{lm}_{\mathbf{d}}(y_{\mathbf{d}\minus1})$: \hspace{8pt} $V^{lm}_\mathbf{d} (y_{\mathbf{d}\minus1}) :=$
\begin{equation}\label{eq:sp_init_D}
	\begin{aligned}
	\underset{y_\mathbf{d} \in \mathcal{Y}_\mathbf{d} }{\text{min}} & \hspace{3pt} \sum_{\omega \in \Omega} \pi^{\omega} \left( c^{\top}_\mathbf{d} y^\omega_\mathbf{d} \right) \\
	\text{s.t.} & \hspace{3pt}  A_\mathbf{d} \hspace{1.0pt} y^\omega_\mathbf{d}  \leq  x^{\top} \hspace{-1pt} B_\mathbf{d} \hspace{1.0pt} b^{l m \omega}_\mathbf{d}  +  C_\mathbf{d} \hspace{1.0pt} y_{\mathbf{d}\minus1}, \;  \omega \in \Omega. \\
	\end{aligned}
\end{equation}

\subsection{Assumptions}\label{subsec:Assump}
In this paper, we assume (1) relatively complete recourse and (2) the uncertainty is stage-wise independent in the subproblems. We assume that problem $\mathrm{SP}^{lm}_d$ defined in~\eqref{eq:sp_init_d} has bounded subgradients, for all $d=1,..,\mathbf{d}$. It follows that there exists a constant $\mathrm{M_y}$ such that $\mathrm{M_y}\geq|| \sigma^{lm}_d\hspace{-1pt}(y_d) ||_1$, for all $y_d \in \mathcal{Y}_d$, for all $l=1,..,\mathbf{m}$, for all $m=1,..,\mathbf{m}$, for all $d=1,..,\mathbf{d}$, where $\sigma^{lm}_d\hspace{-1pt}(y_d)$ is the subgradient of $V^{lm}_d$ w.r.t. $y_d$. Finally, let us denote the norm $||\cdot||_1$ as $||\cdot||$.    

\section{SDDP algorithm with valid bounds and a path search rule}\label{sec:sddp_approx}
In this section, we focus on developing the skeleton of our SDDP algorithm. We initially forget the need to generate valid Benders cut from the solution of our subproblem, and we also ignore the tractability of the proposed algorithm. We also show the search paths that the SDDP algorithm needs to follow in order to prove convergence in a finite number of iterations to an $\delta$-optimal solution.

\subsection{bounding envelopes}\label{subsec:bound_envs}

\begin{lemma}\label{lemma:convy}
	$V^{mn}_{d\plus1}(y^\omega_{d})$ is a convex function of $y^\omega_{d}$, $  d=1,..,\mathbf{d}\minus1$.
\end{lemma}
\begin{proof}
    Let us now assume that $V^{mn}_{d\plus1}(y^\omega_{d})$ is a convex function of $y^\omega_d$. It follows that problem $\mathrm{SP}^{lm}_d$ defined in~\eqref{eq:sp_init_d} is a convex optimisation problem, and that $V^{lm}_d(y_{d\minus1})$ is convex w.r.t. $y_{d\minus1}$ as it only appears as a right-hand side coefficient. Then, observe that $V^{mn}_\mathbf{d}(y^\omega_{\mathbf{d}\minus1})$ is convex w.r.t. $y^\omega_{\mathbf{d}-1}$, which implies that also $V^{lm}_{\mathbf{d}-1}(y_{\mathbf{d}\minus2})$ is convex w.r.t. $y_{\mathbf{d}\minus2}$, and by induction that $V^{mn}_{d\plus1}(y^\omega_{d})$ is a convex function of $y^\omega_{d}$ for all $ d=1,..,\mathbf{d}\minus1$. 
    \qed
\end{proof}

\subsubsection*{Envelope for valid lower bound}

Let us assume we know $\{(\underline{\theta}^{mn}_{s,d\plus1},\underline{\sigma}^{mn}_{s,d\plus1},y^{mn}_{s,d}), s \in \mathcal{S}\}$, such that $\underline{\theta}^{mn}_{s,d\plus1}$ is a valid lower bound on $V^{mn}_{d\plus1}(y^{mn}_{s,d})$, and $\underline{\sigma}^{mn}_{s,d\plus1}$ is a valid subgradient w.r.t. $y^{mn}_{s,d}$.
We define a valid lower bound $\underline{\vartheta}^{mn}_{d\plus1} (y^{\omega}_d)$ on the true value $V^{mn}_{d\plus1} (y^{\omega}_d)$ as
\begin{equation}\label{eq:varthetaLB}
	\begin{aligned}
		\underline{\vartheta}^{mn}_{d\plus1} (y^\omega_d) := \quad \underset{\underline{\vartheta}^{\omega mn}_{d\plus1}}{\text{min}} & \hspace{3pt} \underline{\vartheta}^{\omega mn}_{d\plus1} \\
		\text{s.t.} & \hspace{3pt} \underline{\vartheta}^{\omega mn}_{d\plus1} \geq \underline{\theta}^{mn}_{s,d\plus1} + \underline{\sigma}^{mn\top}_{s,d\plus1} (y^\omega_d - y^{mn}_{s,d}), \hspace{6pt}  s  \in \mathcal{S}. \\
	\end{aligned}
\end{equation}

\begin{lemma}\label{lemma:Venvelope_LB}
	Problem~\eqref{eq:varthetaLB} satisfies the following properties 
	\begin{itemize}
		\item[$i$)] $\underline{\vartheta}^{mn}_{d\plus1}(y^\omega_d) \leq V^{mn}_{d\plus1} (y^\omega_d)$, $ y^\omega_d \in \mathcal{Y}^\omega_d, \omega \in \Omega$,
		\item[$ii$)] $\underline{\vartheta}^{mn}_{d\plus1}(y^\omega_d) + \sigma^{mn}_{d\plus1}(y^\omega_d)^\top (\tilde{y}_d-y^\omega_d) \leq  V^{mn}_{d\plus1} (\tilde{y}_d)$, $ \tilde{y}_d \in \mathcal{Y}^\omega_d$, $ y^\omega_d \in \mathcal{Y}^\omega_d,  \omega \in \Omega$, 
		\item[$iii$)] $\underline{\vartheta}^{mn}_{d\plus1}(y^\omega_d)  \geq \underline{\theta}^{mn}_{s,d\plus1} - \mathrm{M_y}||y^\omega_d-y^{mn}_{s,d}|| $, $ y^\omega_d \in \mathcal{Y}^\omega_d,  \omega \in \Omega,  s \in \mathcal{S}$,
	\end{itemize}
	where $\sigma^{mn}_{d\plus1}(y^\omega_d)$ is the subgradient of $\underline{\vartheta}^{mn}_{d\plus1} (y^\omega_d)$ w.r.t. $y^\omega_d$.
\end{lemma}
\begin{proof} \hspace{1pt} \\
	$i$) Each cut $\underline{\theta}^{mn}_{s,d\plus1} + \underline{\sigma}^{mn\top}_{s,d\plus1} (y^\omega_d - y^{mn}_{s,d})$ of~\eqref{eq:varthetaLB} is lower or equal to $V^{mn}_{d\plus1} (y^\omega_d)$ given that $V^{mn}_{d\plus1} (y^\omega_d)$ is convex w.r.t. $y^\omega_d$. It follows that $\underline{\vartheta}^{mn}_{d\plus1}(y^\omega_d) \leq V^{mn}_{d\plus1} (y^\omega_d)$, $ y^\omega_d \in \mathcal{Y}^\omega_d,  \omega \in \Omega$.\\
	$ii$) $\underline{\vartheta}^{mn}_{d\plus1}(y^\omega_d)$ is convex w.r.t. $y^\omega_d$ since $y^\omega_d$ only appear as a right-hand side coefficient in~\eqref{eq:varthetaLB}. It follows that 
	$$\underline{\vartheta}^{mn}_{d\plus1}(y^\omega_d) + \sigma^{mn}_{d\plus1}(y^\omega_d)^\top (\tilde{y}_d-y^\omega_d) \leq \underline{\vartheta}^{mn}_{d\plus1}(\tilde{y}_d), \quad  \tilde{y}_d \in \mathcal{Y}^\omega_d,  y^\omega_d \in \mathcal{Y}^\omega_d,  \omega \in \Omega$$
	and from part $i$) follows that $\underline{\vartheta}^{mn}_{d\plus1}(\tilde{y}_d)\leq V^{mn}_{d\plus1}(\tilde{y}_d)$, $ \tilde{y}_d \in \mathcal{Y}^\omega_d,  \omega \in \Omega$. \\ 
    $iii$) The definition of $\underline{\vartheta}^{mn}_{d\plus1}(y^\omega_d)$ in~\eqref{eq:varthetaLB}, yields to
	\begin{equation*}\label{eq:UBtheorem_i}
		\begin{aligned}
			\underline{\vartheta}^{mn}_{d\plus1}(y^\omega_d)
			& \geq \hspace{3pt}  \underline{\theta}^{mn}_{s,d\plus1} + \underline{\sigma}^{mn\top}_{s,d\plus1} (y^\omega_d- y_{s,d}) \\
			& \geq \hspace{3pt}  \underline{\theta}^{mn}_{s,d\plus1} - \mathrm{M_y}||y^\omega_d- y^{mn}_{s,d}|| 
		\end{aligned}
	\end{equation*}
    for each $y^\omega_d \in \mathcal{Y}^\omega_d$, for each $\omega \in \Omega$, and for each $s=1,..\mathbf{s}$. The first inequality holds since setting $\underline{\vartheta}^{mn}_{d\plus1}(y^\omega_d)$ is convex w.r.t. $y^\omega_d$, and the second inequality holds since $||\underline{\sigma}^{mn}_{s,d\plus1}||\leq \mathrm{M_y}$.
	\qed 
\end{proof}

A valid lower bound approximation on problem $\mathrm{SP}^{lm}_d(y_{d\minus1})$ of~\eqref{eq:sp_init_d} can then be obtain solving problem $\mathrm{LP}^{lm}_{d}(y_{d\minus1})$ defined as \\

\noindent $\mathrm{LP}^{lm}_{d}(y_{d\minus1}): \hspace{8pt} \underline{\theta}^{lm}_{d} (y_{d\minus1}) :=$
\vspace{-3pt}
\begin{equation}\label{eq:LB_lm}
	\begin{aligned}
	\underset{\substack{y_d \in \mathcal{Y}_d,\\\underline{\vartheta}^{\omega mn}_{d\plus1}}}{\text{min}} & \hspace{3pt} \sum_{\omega \in \Omega} \pi^{\omega} \hspace{-1pt} \left(  c^{\top}_d y^\omega_d + \sum_{n \in \mathcal{M}} \pi^{mn}_{d\plus1} \underline{\vartheta}^{\omega mn}_{d\plus1} \right) \\
	\text{s.t.} & \hspace{3pt} A_d \hspace{1.5pt} y^\omega_{d} \leq x^{\top} \hspace{-1pt} B_d \hspace{1.5pt} b^{lm \omega}_d + C_d \hspace{1.5pt} y_{d\minus1}, \;  \omega \in \Omega \\
                & \hspace{3pt} \underline{\vartheta}^{\omega mn}_{d\plus1} \geq \underline{\theta}^{mn}_{s,d\plus1} + \underline{\sigma}_{s,d\plus1}^{mn\top} (y^\omega_d-y^{mn}_{s,d}), \;  \omega \in \Omega,  n \in \mathcal{M},  s \in \mathcal{S} \\
	\end{aligned}
\end{equation}
Each element $(\underline{\theta}^{mn}_{s,d\plus1},\underline{\sigma}^{mn}_{s,d\plus1},y^{mn}_{s,d})$ is generated solving $\mathrm{LP}^{mn}_{d\plus1}(y^{mn}_{s,d})$ and obtaining the optimal objective $\underline{\theta}^{mn}_{s,d\plus1}$, and a subgradient $\underline{\sigma}_{d\plus1}^{mn}$ w.r.t. $y^{mn}_{s,d}$.
At stage $\mathbf{d}$, $\mathrm{LP}^{lm}_\mathbf{d}(y_{\mathbf{d}\minus1})$ is equivalent to $\mathrm{SP}^{lm}_\mathbf{d}(y_{\mathbf{d}\minus1})$ defined in~\eqref{eq:sp_init_D}.

\begin{lemma}\label{lemma:LP_Lemma}
	Solving problem $\mathrm{LP}^{lm}_{d}(y_{d-1})$ yields an optimal solution $\underline{y}^*_d$ and an optimal objective $\underline{\theta}^*_d$ which is a valid lower bound on the objective of $\mathrm{SP}^{lm}_{d}(y_{d-1})$.
\end{lemma}
\begin{proof} \hspace{1pt} \\
	We start noticing that this is valid for stage $\mathbf{d}$, given that $\mathrm{LP}^{lm}_\mathbf{d}(y_{\mathbf{d}\minus1})$ is equivalent to $\mathrm{SP}^{lm}_\mathbf{d}(y_{\mathbf{d}\minus1})$ by definition. Then, at stage $\mathbf{d}-1$, each element $(\underline{\theta}^{mn}_{s,\mathbf{d}},\underline{\sigma}^{mn}_{s,\mathbf{d}},y^{mn}_{s,\mathbf{d}-1})$ added to $\mathrm{LP}^{lm}_{\mathbf{d}\minus1}$
	is generated solving $\mathrm{LP}^{mn}_\mathbf{d}(y^{mn}_{s,\mathbf{d}-1})$ which is equivalent to $\mathrm{SP}^{mn}_\mathbf{d}(y^{mn}_{s,\mathbf{d}-1})$, and by Lemma~\ref{lemma:Venvelope_LB} follows that the optimal objective $\underline{\theta}^*_{\mathbf{d}\minus1}$ of $\mathrm{LP}^{lm}_{\mathbf{d}\minus1}(y_{\mathbf{d}\minus2})$ is a valid lower bound on $\mathrm{SP}^{lm}_{\mathbf{d}\minus1}(y_{\mathbf{d}\minus2})$. By induction, it follows that $\underline{\theta}^*_d$ is valid lower bound on the objective of $\mathrm{SP}^{lm}_d(y_{d\minus1})$, for all $d=1,..,\mathbf{d}$. \qed
\end{proof}

\subsubsection*{Envelope for valid upper bound}

Let us assume we know $\{(\overline{\theta}^{mn}_{s,d\plus1},y^{mn}_{s,d}), s \in \mathcal{S}\}$, such that $\overline{\theta}^{mn}_{s,d\plus1}$ is a valid upper bound on $V^{mn}_{d\plus1}(y^{mn}_{s,d})$. 
We define a valid upper bound $\overline{\vartheta}^{mn}_{d\plus1} (y^\omega_d)$ on the true value $V^{mn}_{d\plus1} (y^\omega_d)$ as
\begin{equation}\label{eq:varthetaUB}
	\begin{aligned}
	\hspace{5pt} \overline{\vartheta}^{mn}_{d\plus1} (y^\omega_d) :=  \underset{\substack{\overline{\vartheta}^{\omega mn}_{d\plus1},\gamma^{\omega mn}_{d\plus1},\\\mu^{\omega mn}_{s,d\plus1} \geq 0}}{\text{min}} & \hspace{3pt} \overline{\vartheta}^{\omega mn}_{d\plus1} \\
	\text{s.t.} & \hspace{3pt} \overline{\vartheta}^{\omega mn}_{d\plus1} \geq
 \sum\nolimits_{s \in \mathcal{S}} \mu^{\omega mn}_{s,d\plus1} \overline{\theta}^{mn}_{s,d\plus1} + \mathrm{M_y} ||\gamma^{\omega mn}_{d\plus1}||  \\ 
	            & \hspace{3pt} \sum\nolimits_{s \in \mathcal{S}} \mu^{\omega mn}_{s,d\plus1} \hspace{1pt} y^{mn}_{s,d} = y^\omega_d + \gamma^{\omega mn}_{d\plus1}, \\ 
				& \hspace{3pt} \sum\nolimits_{s \in \mathcal{S}} \mu^{\omega mn}_{s,d\plus1} = 1.\\ 
	\end{aligned}
\end{equation}

\begin{lemma}\label{lemma:Venvelope_UB}
	Problem~\eqref{eq:varthetaUB} satisfies the following properties 
	\begin{itemize}
		\item[$i$)]  $\overline{\vartheta}^{mn}_{d\plus1}(y^\omega_d) \geq V^{mn}_{d\plus1} (y^\omega_d)$, $ y^\omega_d \in \mathcal{Y}^\omega_d,  \omega \in \Omega$, 
		\item[$ii$)] $\overline{\vartheta}^{mn}_{d\plus1}(y^\omega_d) \leq \overline{\theta}^{mn}_{s,d\plus1} + \mathrm{M_y}||y^\omega_d-y^{mn}_{s,d}|| $, $ y^\omega_d \in \mathcal{Y}^\omega_d,  \omega \in \Omega , s \in \mathcal{S}$.
	\end{itemize}
\end{lemma}
\begin{proof} \hspace{1pt} \\
	$i$) The definition of $\overline{\vartheta}^{mn}_{d\plus1} (y^\omega_d)$ in~\eqref{eq:varthetaUB} leads to
	\begin{equation*}
		\begin{aligned}
			\overline{\vartheta}^{mn}_{d\plus1} (y^\omega_d)
			& \geq \hspace{3pt} \sum\nolimits_{s \in \mathcal{S}} \mu^{\omega mn}_{s,d\plus1} V^{mn}_{d\plus1} (y^{mn}_{s,d}) \hspace{-1.5pt}+\hspace{-1.5pt} \mathrm{M_y} ||\gamma^{\omega mn}_{d\plus1}|| \\
			& \geq \hspace{3pt} V^{mn}_{d\plus1} (y^\omega_d + \gamma^{\omega mn}_{d\plus1}) + \mathrm{M_y} |||\gamma^{\omega mn}_{d\plus1}||\\
			& \geq \hspace{3pt} V^{mn}_{d\plus1} (y^\omega_d).
		\end{aligned}
	\end{equation*}
	The first inequality holds since $V^{mn}_{d\plus1} (y^{mn}_{s,d}) \leq \overline{\theta}^{mn}_{s,d\plus1}$ for each $s \in \mathcal{S}$, the second inequality holds since $\sum\nolimits_{s \in \mathcal{S}} \mu^{\omega mn}_{s,d\plus1} y^{mn}_{s,d} = y^\omega_d + \gamma^{\omega mn}_{d\plus1}$, the $\mu^{\omega mn}_{s,d\plus1}$ define a convex combination, and $V^{mn}_{d\plus1} (y^\omega_{d})$ is convex w.r.t. $y^\omega_{d}$. The last inequality holds since $||\sigma^{mn}_{s,d\plus1}(y^\omega_d)||\leq \mathrm{M_y}$, $ y^\omega_d \in \mathcal{Y}^\omega_d,  \omega \in \Omega$.\\
	$ii$) The definition of $\overline{\vartheta}^{mn}_{d\plus1} (y^\omega_d)$ in~\eqref{eq:varthetaUB} yields to
	\begin{equation*}
		\begin{aligned}
			\overline{\vartheta}^{mn}_{d\plus1} (y^\omega_d)
			& \leq \hspace{3pt} \overline{\theta}^{mn}_{s,d\plus1} + \mathrm{M_y} ||\gamma^{\omega mn}_{d\plus1}|| \\
			& = \hspace{3pt} \overline{\theta}^{mn}_{s,d\plus1} + \mathrm{M_y}||y^\omega_d-y^{mn}_{s,d}|| 
		\end{aligned}
	\end{equation*}
	The first inequality holds since setting $\mu^{\omega mn}_{s,d\plus1}$ equal to 1 gives a feasible (but not necessarily optimal) solution with objective $\overline{\theta}^{mn}_{s,d\plus1} + \mathrm{M_y} ||\gamma^{\omega mn}_{d\plus1}||$. The equality holds since $\gamma^{\omega mn}_{d\plus1}=y^{mn}_{s,d}-y^\omega_{d}$ when $\mu^{\omega mn}_{s,d\plus1}=1$.
	\qed
\end{proof}

A valid upper bound approximation on problem $\mathrm{SP}^{lm}_d(y_{d\minus1})$ of~\eqref{eq:sp_init_d} can then be obtain solving problem $\mathrm{UP}^{lm}_{d}(y_{d\minus1})$ defined as \\

\noindent $\mathrm{UP}^{lm}_d(y_{d\minus1}): \hspace{8pt} \overline{\theta}^{lm}_{d} (y_{d\minus1}) := $
\vspace{-3pt}
\begin{equation}\label{eq:UBb_1}
	\begin{aligned}
	\underset{\substack{y^d \in \mathcal{Y}_d,\\\overline{\vartheta}^{\omega mn}_{d\plus1},\gamma^{\omega mn}_{d\plus1},\\\mu^{\omega mn}_{s,d\plus1} \geq 0}}{\text{min}} & \hspace{3pt} \sum_{\omega \in \Omega} \pi^\omega \left( c_d^{\top} y^\omega_d + \sum_{n \in \mathcal{M}} \pi^{mn}_{d\plus1} \overline{\vartheta}^{\omega mn}_{d\plus1} \right) \\
	\text{s.t.} & \hspace{3pt} A_d \hspace{1.5pt} y^\omega_{d} \leq x^{\top} \hspace{-1pt} B \hspace{1.5pt} b^{l m \omega}_{d} + C_d \hspace{1.5pt} y_{d\minus1}, \;  \omega \in \Omega,\\
				& \hspace{3pt}\overline{\vartheta}^{\omega mn}_{d\plus1} \geq \sum\nolimits_{s \in \mathcal{S}} \mu^{\omega mn}_{s,d\plus1} \overline{\theta}^{mn}_{s,d\plus1} + \mathrm{M_y} ||\gamma^{\omega mn}_{d\plus1}||, \;  \omega \in \Omega,  n \in \mathcal{M},\\
				& \hspace{3pt} \sum\nolimits_{s \in \mathcal{S}} \mu^{\omega mn}_{s,d\plus1} y^{mn}_{s,d} = y^\omega_d + \gamma^{\omega mn}_{d\plus1}, \;  \omega \in \Omega,  n \in \mathcal{M}, \\ 
				& \hspace{3pt} \sum\nolimits_{s \in \mathcal{S}} \mu^{\omega mn}_{s,d\plus1} = 1, \;  \omega \in \Omega,  n \in \mathcal{M}. 
	\end{aligned}
\end{equation}
Each element $(\overline{\theta}^{mn}_{s,d\plus1},y^{mn}_{s,d})$ is generated solving problem $\mathrm{UP}^{mn}_{d\plus1}(y^{mn}_{s,d})$ and obtaining the optimal objective $\overline{\theta}^{mn}_{s,d\plus1}$. 
At the stage $\mathbf{d}$, problem $\mathrm{UP}^{lm}_\mathbf{d}(y_{\mathbf{d}\minus1})$ is equivalent to $\mathrm{SP}^{lm}_\mathbf{d}(y_{\mathbf{d}\minus1})$ defined in~\eqref{eq:sp_init_D}.
\begin{lemma}\label{lemma:UP_Lemma}
	Solving problem $\mathrm{UP}^{lm}_{d}(y_{d-1})$ yields an optimal solution $\overline{y}^*_d$ and an optimal objective $\overline{\theta}^*_d$ which is a valid upper bound on the objective of $\mathrm{SP}^{lm}_{d}(y_{d-1})$.
\end{lemma}
\begin{proof} \hspace{1pt} \\
	This Lemma can be proved with the same arguments of Lemma~\ref{lemma:LP_Lemma}. \qed
\end{proof}

We define an additional valid bound approximation on $\mathrm{SP}^{lm}_d(y_{d\minus1})$. It is obtained solving problem $\mathrm{UP}^{lm}_{d}(y_{d\minus1})$ where we impose $y^\omega_d = \tilde{y}^\omega$, $ \omega \in \Omega$, i.e., we generate a valid upper bound $\overline{\theta}^{lm}_d$ at point $\tilde{y}$. We define this problem as $\mathrm{UP}^{lm}_{d}(y_{d\minus1},\tilde{y})$.

\subsection{SDDP algorithm}
The SDDP algorithm is presented in Algorithm~\ref{alg:SDDPv0}. In the following, we show the convergence proof of the basic SDDP algorithm.

\begin{algorithm}[!htb]
	\caption{\texttt{SDDP with valid bounds and a path search rule}}
	\label{alg:SDDPv0}
	\algotitle{\texttt{Basic\_SDDP}}{Stand_Bend.title}
	choose tolerance $\delta > 0$ and set $k:=0$\;
	set $\underline{\theta}^k:=-\infty$, $\overline{\theta}^k:=\infty$, and set $\underline{\vartheta}^{m}_1:=-\infty$, $\overline{\vartheta}^{m}_1:=\infty, {\scriptstyle  m \in \mathcal{M}}$\; 
	\Repeat{$\overline{\theta}^{k}-\underline{\theta}^{k}\leq\delta$}{
		set $k:=k+1$, $d:=1$, and $\hat{d}:=\mathbf{d}$\;
		\tcc{Forward Pass}
		set $m_1:=$\raisebox{.75pt}{\scalebox{.90}{$\mathrm{arg} \, \underset{m \in \mathcal{M}}{\mathrm{max}}  \left\{ \left( \overline{\vartheta}^m_1 \hspace{-3pt} - \underline{\vartheta}^m_1 \hspace{-1pt} \right) \,\middle|\, \pi^m_1 > 0 \right\}$}}\;
		solve $\mathrm{LP}^{m_1}_{1}$, get $\underline{y}_1$ and $\{\underline{\vartheta}^{\omega m_1 n}_{2}, {\scriptstyle \omega \in \Omega, n \in \mathcal{M}}\}$\;
		solve $\mathrm{UP}^{m_1}_{1}(\underline{y}_1)$, get $\{\overline{\underline{\vartheta}}^{\omega m_1 n}_{2}, {\scriptstyle \omega \in \Omega, n \in \mathcal{M}}\}$\;
		\Repeat{$\hat{d}=d-1$ $\mathbf{or}$ $d=\mathbf{d}$}{
			set $d:=d+1$\;
			set \raisebox{.75pt}{\scalebox{.90}{$(\omega_{d\minus1},m_d):= \mathrm{arg} \underset{\substack{\omega \in \Omega, n \in \mathcal{M} }}{\mathrm{max}} \left\{ \left( \overline{\underline{\vartheta}}^{\omega m_{d\minus1} n}_d \hspace{-3pt} - \underline{\vartheta}^{\omega m_{d\minus1} n}_d \right) \,\middle|\, \pi^{m_{d\minus1} n}_d > 0 \right\}$}}\;
			\uIf{\raisebox{.75pt}{\scalebox{.90}{$\left( \overline{\underline{\vartheta}}^{\omega_{d\minus1} m_{d\minus1} m_d}_d \hspace{-3pt} - \underline{\vartheta}^{\omega_{d\minus1} m_{d\minus1} m_d }_d \hspace{-1pt} \right)$}}$ \, > \tfrac{\delta (\mathbf{d}-d+1)}{\mathbf{d}-1} $}{
				solve $\mathrm{LP}^{m_{d\minus1} m_d}_{d}(\underline{y}^{\omega_{d\minus1}}_{d\minus1})$, get $\underline{y}_d$ and $\{\underline{\vartheta}^{\omega m_d n}_{d\plus1},{\scriptstyle \omega \in \Omega, \hspace{2pt}  n \in \mathcal{M}}\}$\; 
				solve $\mathrm{UP}^{m_{d\minus1} m_d}_{d}(y^{\omega_{d\minus1}}_{d\minus1},\underline{y}_d)$, get $\{\overline{\underline{\vartheta}}^{\omega m_d n}_{d\plus1},{\scriptstyle \omega \in \Omega, \hspace{2pt}  n \in \mathcal{M}}\}$\;
				}\Else{
				$\hat{d}:=d-1$\;	
				}
		}
		\tcc{Backward Pass}
		\For{$d=\hat{d},\hat{d}-1,..,2$}{
			solve $\mathrm{LP}^{m_{d\minus1} m_{d}}_{d}(\underline{y}^{\omega_{d\minus1}}_{d\minus1})$, get $\underline{\theta}^{m_{d\minus1} m_{d}}_{d}$ and $\underline{\sigma}^{m_{d\minus1} m_{d}}_{d}$\;
			solve $\mathrm{UP}^{m_{d\minus1} m_{d}}_{d}(\underline{y}^{\omega_{d\minus1}}_{d\minus1})$, get $\overline{\theta}^{m_{d\minus1} m_{d}}_{d}$\;
			add cut $(\underline{\theta}^{m_{d\minus1} m_{d}}_{d},\underline{\sigma}^{m_{d\minus1} m_{d}}_{d},\underline{y}^{\omega_{d\minus1}}_{d\minus1})$ to $\mathrm{LP}^{m_{d\minus2} m_{d\minus1}}_{d\minus1}$\;
			add cut $(\overline{\theta}^{m_{d\minus1} m_{d}}_{d},\underline{y}^{\omega_{d\minus1}}_{d\minus1})$ to $\mathrm{UP}^{m_{d\minus2} m_{d\minus1}}_{d\minus1}$\;
		}
        solve $\mathrm{LP}^{m_1}_{1}$ and $\mathrm{UP}^{m_1}_{1}$ to get $\underline{\theta}^{m_{1}}_{1}$ and $\overline{\theta}^{m_1}_{1}$\;
		set $\underline{\vartheta}^{m_1}_1:=\underline{\theta}^{m_{1}}_{1}$ and $\overline{\vartheta}^{m_1}_1:=\overline{\theta}^{m_{1}}_{1}$ \;
		set $\underline{\theta}^{k} := \sum_{m \in \mathcal{M}} \pi^m_1\underline{\vartheta}^{m}_1$, and set $\overline{\theta}^{k}  := \sum_{m \in \mathcal{M}}\pi^m_1\overline{\vartheta}^{m}_1$ \;
	}
\end{algorithm}

\subsubsection*{Convergence proof of the basic SDDP algorithm}
\definition\label{def:suset_R}
Let $\{(\underline{\theta}^{mn}_{s,d\plus1},\underline{\sigma}^{mn}_{s,d\plus1},y^{mn}_{s,d}),s \in \mathcal{S}^{mn}_{d\plus1}\}$ and $\{(\overline{\theta}^{mn}_{s,d\plus1},y^{mn}_{s,d}),s \in \mathcal{S}^{mn}_{d\plus1}\}$ be the sets of elements already added to $\mathrm{LP}^{lm}_d$ and $\mathrm{UP}^{lm}_d$. We define $R^{mn}_{d\plus1}$ as the subset of $S^{mn}_{d\plus1}$, that satisfies
\begin{itemize}
	\item[$i)$] $\overline{\theta}^{mn}_{s,d\plus1}-\underline{\theta}^{mn}_{s,d\plus1} \leq \frac{\delta(\mathbf{d}-d-1)}{\mathbf{d}-1}$, $ s \in R^{mn}_{d\plus1}$,  
	\item[$ii)$] $||y^{mn}_{s,d}-y^{mn}_{r,d}|| > \frac{\delta}{2(\mathbf{d}-1)\mathrm{M_y}}$, $ s,r \in R^{mn}_{d\plus1}$ such that $s \neq r$,    
\end{itemize}
for each $n \in \mathcal{M}$ and for each $d=1,..,\mathbf{d}-1$.

\begin{lemma}\label{lemma:subset_R_finite}
	There exists a finite number of elements that can be added to $R^{mn}_{d\plus1}$, for each $n \in \mathcal{M}$ and for each $d=1,..,\mathbf{d}-1$.
\end{lemma}
\begin{proof} \hspace{1pt} \\
	Since $\delta >0$, $\mathrm{M_y}$ is finite, and $\mathbf{d}$ is finite, it follows that $0 <\frac{\delta}{2(\mathbf{d}-1)\mathrm{M_y}}<\infty$. Given that each $y^{mn}_{s,d}$ belongs to a compact set, it follows that there exists a finite number of elements that can be added to $R^{mn}_{d\plus1}$, for each $n \in \mathcal{M}$ and for each $d=1,..,\mathbf{d}-1$. \qed
\end{proof}

\begin{lemma}\label{lemma:accuracy_backward}
	Let $\underline{y}_{d}$ and $\{\underline{\vartheta}^{\omega mn}_{d\plus1}, \omega \in \Omega, n \in \mathcal{M} \}$ be the optimal solution of problem $\mathrm{LP}^{lm}_{d}(y_{d\minus1})$, and $\underline{\theta}^{lm}_{d}$ its optimal objective. Then, let $\{\overline{\vartheta}^{\omega mn}_{d\plus1}, \omega \in \Omega, n \in \mathcal{M} \}$ be the optimal solution of problem $\mathrm{UP}^{lm}_{d}(y_{d\minus1})$, and $\overline{\theta}^{lm}_{d}$ its optimal objective.
	Finally, let $\{\overline{\underline{\vartheta}}^{\omega mn}_{d\plus1}, \omega \in \Omega, n \in \mathcal{M} \}$ be the optimal solution of problem $\mathrm{UP}^{lm}_{d}(y_{d\minus1},\underline{y}_{d})$. If
	$$\hspace{-2.5pt} \left(\overline{\underline{\vartheta}}^{\tilde{\omega}m\tilde{n}}_{d\plus1} - \underline{\vartheta}^{\tilde{\omega}m\tilde{n}}_{d\plus1}\right) \leq \tfrac{\delta (\mathbf{d}-d)}{ \mathbf{d}-1 },$$
	where 
	$$ (\tilde{\omega},\tilde{n}) :=  \mathrm{arg} \underset{\substack{ \omega \in \Omega, n \in \mathcal{M} }}{\mathrm{max}} \left(  \overline{\underline{\vartheta}}^{\omega mn}_{d\plus1} \hspace{-3pt} - \underline{\vartheta}^{\omega mn}_{d\plus1} \right) \; | \; \pi^{mn}_d > 0 $$
	it follows that $\overline{\theta}^{lm}_{d}-\underline{\theta}^{lm}_{d}\leq \tfrac{\delta (\mathbf{d}-d)}{(\mathbf{d}-1)}$.
\end{lemma}
\begin{proof} \hspace{1pt} \\
	The definition of $\overline{\theta}^{lm}_{d}$ and $\underline{\theta}^{lm}_{d}$ leads to

	\begin{equation*}
        \begin{aligned}
            \overline{\theta}^{lm}_{d} - \underline{\theta}^{lm}_{d}
			= & \sum_{\omega \in \Omega} \pi^{\omega} \hspace{-1pt} \left( \left(  c_d^\top \underline{y}^\omega_d + \hspace{-3pt} \sum_{n \in \mathcal{M}} \pi^{mn}_{d\plus1} \overline{\vartheta}^{\omega mn}_{d\plus1} \right) - \left(c_d^\top \underline{y}^\omega_d + \hspace{-3pt} \sum_{n \in \mathcal{M}} \pi^{mn}_{d\plus1} \underline{\vartheta}^{\omega mn}_{d\plus1} \right) \right) \\
			= & \sum_{\omega \in \Omega} \pi^{\omega} \sum_{n \in \mathcal{M}} \pi^{mn}_d \left(\overline{\vartheta}^{\omega mn}_{d\plus1} - \underline{\vartheta}^{\omega mn}_{d\plus1}\right) \\
			\leq & \sum_{\omega \in \Omega} \pi^{\omega} \sum_{n \in \mathcal{M}} \pi^{mn}_d \left(\overline{\underline{\vartheta}}^{\omega mn}_{d\plus1} - \underline{\vartheta}^{\omega mn}_{d\plus1}\right) \\
			\leq & \left(\overline{\underline{\vartheta}}^{\tilde{\omega}m\tilde{n}}_{d\plus1} - \underline{\vartheta}^{\tilde{\omega}m\tilde{n}}_{d\plus1}\right)  \\
			\leq & \; \tfrac{\delta (\mathbf{d}-d)}{(\mathbf{d}-1)} \\
        \end{aligned}
    \end{equation*}

	The first inequality holds since $\{\overline{\underline{\vartheta}}^{\omega mn}_{d\plus1}, \omega \in \Omega, n \in \mathcal{M} \}$ is a feasible but not necessarily optimal solution of problem $\mathrm{UP}^{lm}_{d}(y_{d\minus1})$, and the second inequality follows from the definition of $(\tilde{n},\tilde{\omega})$ and from $\sum_{n\in \mathcal{M}} \pi^{mn}_d$ and $\sum_{\omega\in\Omega} \pi^{\omega}$ being both equal to 1. \qed
\end{proof}

\begin{lemma}\label{lemma:accuracy_R}
	Let $\underline{y}_{d}$ and $\{\underline{\vartheta}^{\omega mn}_{d\plus1}, \omega \in \Omega, n \in \mathcal{M} \}$ be the optimal solution of problem $\mathrm{LP}^{lm}_{d}(y_{d\minus1})$, and let $\{\overline{\underline{\vartheta}}^{\omega mn}_{d\plus1}, \omega \in \Omega, n \in \mathcal{M} \}$ be the optimal solution of problem $\mathrm{UP}^{lm}_{d}(y_{d\minus1},\underline{y}_{d})$. 
	If it exists $r \in R^{mn}_{d\plus1}$ such that $ ||\underline{y}^\omega_{d}-y^{mn}_{r,d}|| \leq \frac{\delta}{2(\mathbf{d}-1)\mathrm{M_y}}$, it follows that
	$$ \overline{\underline{\vartheta}}^{\omega mn}_{d\plus1} - \underline{\vartheta}^{\omega mn}_{d\plus1} \leq \tfrac{\delta (\mathbf{d}-d)}{\mathbf{d}-1}.$$
\end{lemma}
\begin{proof} \hspace{1pt} \\
	We write
	\begin{equation*}
        \begin{aligned}
            \overline{\underline{\vartheta}}^{\omega mn}_{d\plus1} - \underline{\vartheta}^{\omega mn}_{d\plus1}
			\leq & \; \left( \overline{\theta}^{mn}_{r,d\plus1} + \mathrm{M_y} ||\underline{y}^\omega_{d}-y^{mn}_{r,d}|| \right) - \left( \underline{\theta}^{mn}_{r_n,d\plus1} - \mathrm{M_y} ||\underline{y}^\omega_{d}-y^{mn}_{r_n,d}|| \right) \\
			\leq & \; \tfrac{\delta(\mathbf{d}-d-1)}{\mathbf{d}-1} + 2 \mathrm{M_y} \tfrac{\delta}{2(\mathbf{d}-1)\mathrm{M_y}} \\
			   = & \; \tfrac{\delta(\mathbf{d}-d)}{\mathbf{d}-1}.  \\
        \end{aligned}
    \end{equation*}
	The first inqeuality follows from part $iii)$ of Lemma~\ref{lemma:Venvelope_LB} and from part $ii)$ of Lemma~\ref{lemma:Venvelope_UB}
\end{proof}

\begin{lemma}\label{lemma:forward_term_hatd}
    At iteration $\hat{k}$, if Algorithm~\ref{alg:SDDPv0} terminates the forward pass for $\hat{d}>1$, it follows that a new element is added to 
    $\mathcal{R}^{m_{\hat{d}\minus1} m_{\hat{d}}}_{\hat{d}}$ in the backward pass.
\end{lemma}

\begin{proof} \hspace{1pt} \\
	If the forward pass of Algorithm~\ref{alg:SDDPv0} stops at $\hat{d}$, it follows that the optimal solution $\underline{y}_{\hat{d}}$ and $\{\underline{\vartheta}^{\omega m_{\hat{d}} n}_{\hat{d}\plus1}, \omega \in \Omega, n \in \mathcal{M} \}$ of $\mathrm{LP}^{m_{\hat{d}\minus1}m_{\hat{d}}}_{\hat{d}}(\underline{y}^{\omega_{\hat{d}\minus1}}_{\hat{d}\minus1})$ and the optimal solution $\{\overline{\underline{\vartheta}}^{\omega m_{\hat{d}} n}_{\hat{d}\plus1}, \omega \in \Omega, n \in \mathcal{M} \}$ of $\mathrm{UP}^{m_{\hat{d}\minus1}m_{\hat{d}}}_{\hat{d}}(\underline{y}^{\omega_{\hat{d}\minus1}}_{\hat{d}\minus1},\underline{y}_{\hat{d}})$ satisfy
	$$   \overline{\underline{\vartheta}}^{\omega_{\hat{d}} m_{\hat{d}} m_{\hat{d}\plus1}}_{\hat{d}\plus1} - \underline{\vartheta}^{\omega_{\hat{d}} m_{\hat{d}} m_{\hat{d}\plus1}}_{\hat{d}\plus1}  \leq \tfrac{\delta (\mathbf{d}-\hat{d})}{\mathbf{d}-1}.$$
	By Lemma~\ref{lemma:accuracy_backward} it follows that solving $\mathrm{LP}^{m_{\hat{d}\minus1}m_{\hat{d}}}_{\hat{d}}(\underline{y}^{\omega_{\hat{d}\minus1}}_{\hat{d}\minus1})$ and $\mathrm{UP}^{m_{\hat{d}\minus1}m_{\hat{d}}}_{\hat{d}}(\underline{y}^{\omega_{\hat{d}\minus1}}_{\hat{d}\minus1})$ gives optimal objectives such that $\overline{\theta}^{m_{\hat{d}\minus1}m_{\hat{d}}}_{\hat{d}}-\underline{\theta}^{m_{\hat{d}\minus1}m_{\hat{d}}}_{\hat{d}}\leq \tfrac{\delta (\mathbf{d}-\hat{d})}{(\mathbf{d}-1)}$. Note that this is also valid for $\hat{d}=\mathbf{d}$, given that $\overline{\vartheta}^{\omega m_{\mathbf{d}} n}_{\mathbf{d}\plus1}=\underline{\vartheta}^{\omega m_{\mathbf{d}} n}_{\mathbf{d}\plus1}=0,  \omega \in \Omega,  n \in \mathcal{M}$. 
	
	The termination of the forward pass at stage $\hat{d}>1$ also implies that the optimal solution $\underline{y}_{\hat{d}\minus1}$ and $\{\underline{\vartheta}^{\omega m_{\hat{d}\minus1} n}_{\hat{d}}, \omega \in \Omega, n \in \mathcal{M} \}$ of $\mathrm{LP}^{m_{\hat{d}\minus2}m_{\hat{d}\minus1}}_{\hat{d}\minus1}(\underline{y}^{\omega_{\hat{d}\minus2}}_{\hat{d}\minus2})$ and the optimal solution $\{\overline{\underline{\vartheta}}^{\omega m_{\hat{d}\minus1} n}_{\hat{d}}, \omega \in \Omega, n \in \mathcal{M} \}$ of $\mathrm{UP}^{m_{\hat{d}\minus2}m_{\hat{d}\minus1}}_{\hat{d}}(\underline{y}^{\omega_{\hat{d}\minus2}}_{\hat{d}\minus2},\underline{y}_{\hat{d}\minus1})$ are such that
	$$ \ \overline{\underline{\vartheta}}^{\omega_{\hat{d}\minus1} m_{\hat{d}\minus1} m_{\hat{d}}}_{\hat{d}} - \underline{\vartheta}^{\omega_{\hat{d}\minus1} m_{\hat{d}\minus1} m_{\hat{d}}}_{\hat{d}} > \tfrac{\delta (\mathbf{d}-\hat{d}+1)}{\mathbf{d}-1}.$$
	By Lemma~\ref{lemma:accuracy_R} it follows that $\left.\middle|\middle| \underline{y}^{\omega_{\hat{d}\minus1}}_{\hat{d}\minus1} -y^{m_{\hat{d}\minus1} m_{\hat{d}}}_{s,\hat{d}\minus1} \middle|\middle|\right. > \frac{\delta}{2(\mathbf{d}-1)\mathrm{M_y}}, \;  s \in \mathcal{R}^{m_{\hat{d}\minus1} m_{\hat{d}}}_{\hat{d}}$. Hence, the cuts $(\underline{\theta}^{m_{\hat{d}\minus1} m_{\hat{d}}}_{\hat{d}},\underline{\sigma}^{m_{\hat{d}\minus1} m_{\hat{d}}}_{\hat{d}},\underline{y}^{\omega_{\hat{d}\minus1}}_{\hat{d}\minus1})$ and 
	$(\overline{\theta}^{m_{\hat{d}\minus1} m_{\hat{d}}}_{\hat{d}},\underline{y}^{\omega_{\hat{d}\minus1}}_{\hat{d}\minus1})$ added in the backward pass to $\mathrm{LP}^{m_{\hat{d}\minus2} m_{\hat{d}\minus1}}_{\hat{d}\minus1}$ and $\mathrm{UP}^{m_{\hat{d}\minus2} m_{\hat{d}\minus1}}_{\hat{d}\minus1}$, respectively, satisfy both conditions of Definition~\ref{def:suset_R} and add a new element to $\mathcal{R}^{m_{\hat{d}\minus1} m_{\hat{d}}}_{\hat{d}}$. \qed
\end{proof}

\begin{lemma}\label{lemma:forward_term_1}
    At iteration $\hat{k}$, if Algorithm~\ref{alg:SDDPv0} terminates the forward pass for $\hat{d}=1$, it follows that $\overline{\vartheta}^{m_1}_1-\underline{\vartheta}^{m_1}_1 < \delta, \;  k \geq \hat{k}$. 
\end{lemma}

\begin{proof} \hspace{1pt} \\
	If the forward pass of Algorithm~\ref{alg:SDDPv0} stops at $\hat{d}=1$, it follows that the optimal solution $\underline{y}_1$ and $\{\underline{\vartheta}^{\omega m_1 n}_{2}, \omega \in \Omega, n \in \mathcal{M} \}$ of $\mathrm{LP}^{m_1}_{1}$ and the optimal solution $\{\overline{\underline{\vartheta}}^{\omega m_{1} n}_{\hat{d}\plus1}, \omega \in \Omega, n \in \mathcal{M}\}$ of $\mathrm{UP}^{m_1}_{1}(\underline{y}_1)$ satisfy
	$$ \overline{\underline{\vartheta}}^{\omega_1 m_1 m_2}_{2} - \underline{\vartheta}^{\omega_1 m_1 m_2}_{2}  \leq \delta.$$
	By Lemma~\ref{lemma:accuracy_backward} it follows that solving $\mathrm{LP}^{m_1}_{1}$ and $\mathrm{UP}^{m_1}_{1}$ gives optimal objectives such that $\overline{\theta}^{m_1}_{1}-\underline{\theta}^{m_1}_{1}\leq \delta$. Hence, $\overline{\vartheta}^{m_1}_1-\underline{\vartheta}^{m_1}_1 < \delta, \;  k \geq \hat{k}$. \qed
\end{proof}

\begin{theorem}\label{theorem:conv_sddp1}
	For given convergence tolerance $\delta>0$, Algorithm~\ref{alg:SDDPv0} converges to an $\delta$-optimal solution in a finite number of iterations.
\end{theorem}
\begin{proof} \hspace{1pt} \\
	By Lemma~\ref{lemma:forward_term_hatd} Algorithm~\ref{alg:SDDPv0} adds a new element to $\mathcal{R}^{m_{\hat{d}\minus1} m_{\hat{d}}}_{\hat{d}}$ in each iteration for which the forward pass stops at $\hat{d}>1$, and by Lemma~\ref{lemma:subset_R_finite} there exists a finite number of elements that can be added to each $\mathcal{R}^{mn}_{d\plus1}$, for $m \in \mathcal{M}$, $n \in \mathcal{M}$, $d=1,..,\mathbf{d}-1$. Given that the number $\mathbf{d}$ of stages and the number $\mathbf{m}$ of stages is finite, for each initial trajectory $m_1 \in \mathcal{M}$ there exists a finite number of iterations $\hat{k}_1$ for which the forward pass of Algorithm~\ref{alg:SDDPv0} can stop at $\hat{d}>1$. Then, the forward pass will stop at $\hat{d}=1$ and by Lemma~\ref{lemma:forward_term_1} follows that $\overline{\vartheta}^{m_1}_1-\underline{\vartheta}^{m_1}_1 < \delta, \;  k \geq \hat{k}_1$. Given the amount of initial trajectories $\mathbf{m}$ is finite, there exists a finite number of iterations $\mathbf{k}$ such that  $\overline{\vartheta}^{m}_1-\underline{\vartheta}^{m}_1 < \delta, \; m \in \mathcal{M},  k \geq \mathbf{k}$. It follows that $\overline{\theta}^{\mathbf{k}}-\underline{\theta}^{\mathbf{k}} \leq \delta$ and Algorithm~\ref{alg:SDDPv0} finds an $\delta$-optimal solution. \qed
\end{proof}

\section{Enhanced SDDP algorithm}\label{sec:sddp_extended}
	This section extends the formulation of Section~\ref{sec:sddp_approx}. First, we treat $x$ as a variable in the lower bounding envelopes, even if its value is fixed to the one imposed by the master problem. This allows the generation of a valid sensitivity $\underline{\lambda}(x)$ w.r.t. $x$ associated with a valid lower bound $\underline{\theta}(x)$ on the optimal solution of $\mathbf{SP}(x)$ in~\eqref{eq:multistage_g}. Then, we treat the random realisation $b^{lm\omega}_d$ in problem $\mathrm{SP}^{lm}_d(y_{d\minus1})$ defined in~\eqref{eq:sp_init_d} as a variable whose value is fixed to the specific realisation. This reformulation allows to reduce the number of subproblems from $(\mathbf{d}-1)\times\mathbf{m}^2+\mathbf{m}$ (when formulated as \eqref{eq:sp_init_0}-\eqref{eq:sp_init_D}) to $\mathbf{d}\times\mathbf{m}$ and it also gives the option of sharing cuts w.r.t. $b^{mn}_{d\plus1} = \{ b^{mn\omega}_{d\plus1}, \omega \in \Omega\}$ among different subproblems of stage $d$.
	
	\subsection{bounding envelopes}\label{subsec:extended_bnd_envelopes}

	Let define $\mathcal{V}^m_1(x)$ as the optimal solution of\\
	
	\noindent $\mathrm{SP}^{m}_{1}(x)$: \hspace{8pt} $\mathcal{V}^{m}_1 (x) :=$
	\begin{equation}\label{eq:sp_xyb_0} 
		\begin{aligned}
			\underset{y_1 \in \mathcal{Y}_1}{\text{min}} & \; \sum_{\omega \in \Omega} \pi^\omega \left( c^{\top}_1 y^{\omega}_1  + \sum_{n\in \mathcal{M}} \hspace{-2pt} \pi^{mn}_{2} \mathcal{V}^{n}_{2} \hspace{-2pt}\left(x,y^{\omega}_1,b^{mn}_{2}\right)\right)\\
			\text{s.t.} & \; A_1  \hspace{1pt} y^{\omega}_{1}  \leq  x^{\top} \hspace{-1pt} B_1 \hspace{1pt} b^{m\omega}_1, \;  \omega \in \Omega \\
		\end{aligned}
	\end{equation}
	where $\mathcal{V}^m_{d}\hspace{-2pt}\left(x,y_{d\minus1},b^{lm}_{d}\right)$ is given by\\

	\noindent $\mathrm{SP}^{m}_{d}(x,y_{d\minus1},b^{lm}_d)$: \hspace{8pt} $\mathcal{V}^{m}_d (x,y_{d\minus1},b^{lm}_d) :=$
	\begin{equation}\label{eq:sp_xyb_d} 
		\begin{aligned}
			\underset{y_d \in \mathcal{Y}_d}{\text{min}} & \; \sum_{\omega \in \Omega} \pi^\omega \left( c^{\top}_d y^\omega_d  + \sum_{n \in \mathcal{M}} \pi^{mn}_{d\plus1} \mathcal{V}^{n}_{d\plus1} \left(x,y^{\omega}_{d},b^{mn}_{d\plus1}\right) \right)\\
			\text{s.t.} & \; A_d  \hspace{1pt} y^\omega_{d}  \leq  x^{\top} \hspace{-1pt} B_d \hspace{1pt} b^{l m \omega}_d  + C_d \hspace{1pt} y_{d\minus1}, \;  \omega \in \Omega \\
		\end{aligned}
	\end{equation}
	and $\mathcal{V}^{m}_{\mathbf{d}}\hspace{-2pt}\left(x,y_{\mathbf{d}\minus1},b^{lm}_{\mathbf{d}}\right)$ is defined as\\
	
	\noindent $\mathrm{SP}^{m}_{\mathbf{d}}(x,y_{\mathbf{d}\minus1},b^{lm}_\mathbf{d})$: \hspace{8pt} $\mathcal{V}^{m}_\mathbf{d} (x,y_{\mathbf{d}\minus1},b^{lm}_\mathbf{d}) :=$
	\begin{equation}\label{eq:sp_xyb_D} 
		\begin{aligned}
			\underset{y_\mathbf{d} \in \mathcal{Y}_\mathbf{d}}{\text{min}} & \; \sum_{\omega \in \Omega} \pi^\omega \left( c^{\top}_\mathbf{d} y^\omega_\mathbf{d} \right) \\
			\text{s.t.} & \; A_\mathbf{d}  \hspace{1pt} y^\omega_{\mathbf{d}}  \leq  x^{\top} \hspace{-1pt} B_\mathbf{d} \hspace{1pt} b^{lm\omega}_\mathbf{d}  + C_\mathbf{d} \hspace{1pt} y_{\mathbf{d}\minus1}, \;  \omega \in \Omega \\
		\end{aligned}
	\end{equation}

	Note that $\mathcal{V}^m_1(x)$ is equivalent to $V^m_1$ defined in~\eqref{eq:sp_init_0}, $\mathcal{V}^{m}_d \hspace{-2pt}\left(x,y_{d\minus1},b^{lm}_{d}\right)$ is equivalent to $V^{lm}_d \hspace{-2pt}\left(y_{d\minus1}\right)$ defined in~\eqref{eq:sp_init_d}, and $\mathcal{V}^{m}_{\mathbf{d}}\hspace{-2pt}\left(x,y_{\mathbf{d}\minus1},b^{lm}_{\mathbf{d}}\right)$ is equivalent to $V^{lm}_{\mathbf{d}}\hspace{-2pt}\left(y_{\mathbf{d}\minus1}\right)$ defined in~\eqref{eq:sp_init_D}

	\begin{lemma}\label{lemma:conVy}
		$\mathcal{V}^{n}_{d\plus1}\hspace{-2pt}\left(x,y_{d},b^{mn}_{d\plus1}\right)$ is a convex function of $y_d$, $  d=1,..,\mathbf{d}-1$.
	\end{lemma}
	\begin{proof}
		This lemma can be proved with the same argument of Lemma~\ref{lemma:convy}. \qed
	\end{proof}
	\begin{lemma}\label{lemma:conVx}
		$\mathcal{V}^{n}_{d\plus1}\hspace{-2pt}\left(x,y_{d},b^{mn}_{d\plus1}\right)$ is a convex function of $x$, $  d=1,..,\mathbf{d}-1$.
	\end{lemma}
	\begin{proof}
		Let us assume that $\mathcal{V}^{n}_{d\plus1}\hspace{-2pt}\left(x,y_{d},b^{mn}_{d\plus1}\right)$ is a convex function of $x$. It follows that problem $\mathbf{SP}^{m}_d$ defined in~\eqref{eq:sp_xyb_d} is a convex optimisation problem given that $\mathcal{V}^{n}_{d\plus1}\hspace{-2pt}\left(x,y_{d},b^{mn}_{d\plus1}\right)$ is also convex w.r.t. $y_d$ by Lemma~\ref{lemma:conVy}. Then, we notice that also $\mathcal{V}^{m}_{d}\hspace{-2pt}\left(x,y_{d\minus1},b^{ln}_{d}\right)$ is convex w.r.t. $x$ as it appears as a right-hand side coefficient and in the objective in the term $\mathcal{V}^{n}_{d\plus1}\hspace{-2pt}\left(x,y_{d},b^{mn}_{d\plus1}\right)$. Finally, we observe that $\mathcal{V}^{n}_\mathbf{d}\hspace{-2pt}\left(x,y_{\mathbf{d}\minus1},b^{mn}_{\mathbf{d}}\right)$ in~\eqref{eq:sp_xyb_D} is convex w.r.t. $x$, which implies that also $\mathcal{V}^{m}_{\mathbf{d}-1}\hspace{-2pt}\left(x,y_{\mathbf{d}\minus2},b^{lm}_{\mathbf{d}\minus1}\right)$ is convex w.r.t. $x$, and by induction that $\mathcal{V}^{n}_{d\plus1}\hspace{-2pt}\left(x,y_{d},b^{mn}_{d\plus1}\right)$ is a convex function of $x$ for all $ d=1,..,\mathbf{d}-1$. \qed
	\end{proof}

	\begin{lemma}\label{lemma:conVb}
		$\mathcal{V}^{n}_{d\plus1}\hspace{-2pt}\left(x,y_{d},b^{mn}_{d\plus1}\right)$ is a convex function of $b^{mn}_{d\plus1}$, $  d=1,..,\mathbf{d}-1$.
	\end{lemma}
	\begin{proof}
		Let us assume that $\mathcal{V}^{n}_{d\plus1}\hspace{-2pt}\left(x,y_{d},b^{mn}_{d\plus1}\right)$ is a convex function of $b^{mn}_{d\plus1}$. It follows that problem $\mathbf{SP}^{m}_d$ defined in~\eqref{eq:sp_xyb_d} is a convex optimisation problem as $\mathcal{V}^{n}_{d\plus1}\hspace{-2pt}\left(x,y_{d},b^{mn}_{d\plus1}\right)$ is also convex w.r.t. $y_d$ by Lemma~\ref{lemma:conVy}, and that $\mathcal{V}^{m}_{d}\hspace{-2pt}\left(x,y_{d\minus1},b^{lm}_{d}\right)$ is convex w.r.t. $b^{lm}_{d}$ as each $b^{lm\omega}_{d}$ only appears as a right-hand side coefficient. Then, we observe that $\mathcal{V}^{n}_\mathbf{d}\hspace{-2pt}\left(x,y_{\mathbf{d}\minus1},b^{mn}_{\mathbf{d}}\right)$ in~\eqref{eq:sp_xyb_D} is convex w.r.t. $b^{mn}_{\mathbf{d}}$, which implies that also $\mathcal{V}^{m}_{\mathbf{d}-1}\hspace{-2pt}\left(x,y_{\mathbf{d}\minus2},b^{lm}_{\mathbf{d}\minus1}\right)$ is convex w.r.t. $b^{lm}_{\mathbf{d}\minus1}$, and by induction that $\mathcal{V}^{n}_{d\plus1}\hspace{-2pt}\left(x,y_{d},b^{mn}_{d\plus1}\right)$ is a convex function of $b^{mn}_{d\plus1}$ for all $ d=1,..,\mathbf{d}-1$. \qed
	\end{proof}

	\subsubsection*{Envelope for valid lower bound}

	Let us assume we know $\{(\underline{\theta}^{n}_{s,d\plus1},\underline{\lambda}^{n}_{s,d\plus1},\underline{\sigma}^{n}_{s,d\plus1},\underline{\nu}^{n}_{s,d\plus1},x^n_s,y^{n}_{s,d},b^{n}_{s,d\plus1}), s \in \mathcal{S}\}$, such that $\underline{\theta}^{n}_{s,d\plus1}$ is a valid lower bound on $\mathcal{V}^{n}_{d\plus1} \big(x^n_s,y^{n}_{s,d},b^{n}_{s,d\plus1}\big)$, $\underline{\lambda}^{n}_{s,d\plus1}$ is a valid subgradient w.r.t. $x^n_s$, $\underline{\sigma}^{n}_{s,d\plus1}$ is a valid subgradient w.r.t. $y^{n}_{s,d}$, and $\underline{\nu}^{n}_{s,d\plus1}$ is a valid subgradient w.r.t. $b^{n}_{s,d\plus1}$.
	A valid lower bound approximation on problem $\mathbf{SP}^{m}_{d}\hspace{-2pt}\left(x,y_{d\minus1},b^{lm}_{d}\right)$ of~\eqref{eq:sp_xyb_d} can then be obtain solving problem $\mathbf{LP}^{m}_{d}(x,y_{d\minus1},b^{lm}_d)$ defined as \\
	
	\noindent $\mathbf{LP}^{m}_{d}\hspace{-2pt}\left(x,y_{d\minus1},b^{lm}_d\right): \hspace{8pt} \underline{\theta}^{m}_{d} \hspace{-2pt}\left(x,y_{d\minus1},b^{lm}_d\right) :=$
	\vspace{-3pt}
	\begin{equation}\label{eq:LB_lm_1}
		\begin{aligned}
		\underset{\substack{y_d \in \mathcal{Y}_d,\\\underline{\vartheta}^{\omega n}_{d\plus1}}}{\text{min}} & \; \sum_{\omega \in \Omega} \pi^\omega \left( c^{\top}_d y^\omega_d + \sum_{n \in \mathcal{M}} \pi^{mn}_{d\plus1} \underline{\vartheta}^{\omega n}_{d\plus1} \right) \\
		\text{s.t.} & \; A_d \hspace{1.5pt} y^\omega_{d} \leq x^{\top} \hspace{-1pt} B_d \hspace{1.5pt} b^{\ell m \omega}_d + C_d \hspace{1.5pt} y_{d\minus1} , \;  \omega \in \Omega \\
					& \; \underline{\vartheta}^{\omega n}_{d\plus1} \geq \underline{\theta}^{n}_{s,d\plus1} + \underline{\lambda}_{s,d\plus1}^{n\top} (x-x_{s}) + \underline{\sigma}_{s,d\plus1}^{n\top} (y^\omega_d-y^{n}_{s,d}) +  \\
					& \hspace{60pt} + \underline{\nu}_{s,d\plus1}^{n\top} (b^{mn}_{d\plus1}-b^{n}_{s,d\plus1}), \;  \omega \in \Omega,  s \in \mathcal{S}, n \in \mathcal{M}. \\
		\end{aligned}
	\end{equation}

	\begin{lemma}\label{lemma:LPbf_Lemma}
		Solving problem $\mathbf{LP}^{m}_{d}(x,y_{d-1},b^{lm}_d)$ yields an optimal solution $\underline{y}_d$ and an optimal objective $\underline{\theta}_d$ which is a valid lower bound on the objective of $\mathbf{SP}^{m}_{d}(x,y_{d-1},b^{lm}_d)$.
	\end{lemma}
	\begin{proof}
		This lemma can be proved with the same argument of Lemma~\ref{lemma:LP_Lemma}. \qed
	\end{proof}

	\subsubsection*{Envelope for valid upper bound}

	Let us assume we know $\{(\overline{\theta}^{n}_{s,d\plus1},y^{n}_{s,d},b^{n}_{s,d\plus1}), s \in \mathcal{S}\}$, such that $\overline{\theta}^{n}_{s,d\plus1}$ is a valid lower bound on $V^{n}_{d\plus1} \big(x,y^{n}_{s,d},b^{n}_{s,d\plus1}\big)$. A valid upper bound approximation on problem $\mathbf{SP}^{m}_{d}\hspace{-2pt}\left(x,y_{d\minus1},b^{lm}_{d}\right)$ of~\eqref{eq:sp_xyb_d} can then be obtain solving problem $\mathbf{UP}^{m}_{d}\hspace{-2pt}\left(y_{d\minus1},b^{lm}_d\right)$ defined as \\

	\noindent $\mathbf{UP}^{m}_d\hspace{-2pt}\left(y_{d\minus1},b^{lm}_d\right): \hspace{8pt} \overline{\theta}^{m}_{d}\hspace{-2pt}\left(y_{d\minus1},b^{lm}_d\right) := $
	\vspace{-3pt}
	\begin{equation}\label{eq:UBb}
		\begin{aligned}
		\underset{\substack{y^d \in \mathcal{Y}_d,\\\overline{\vartheta}^{\omega n}_{d\plus1},\gamma^{\omega n}_{d\plus1},\\\mu^{\omega n}_{s,d\plus1} \geq 0}}{\text{min}} & \;  \sum_{\omega \in \Omega} \pi^\omega \left( c^{\top}_d y^\omega_d + \sum_{n \in \mathcal{M}} \pi^{mn}_{d\plus1} \overline{\vartheta}^{\omega n}_{d\plus1} \right) \\
		\text{s.t.} & \; A_d \hspace{1.5pt} y^\omega_{d} \leq x^{\top} \hspace{-1pt} B \hspace{1.5pt} b^{l m \omega}_{d} + C_d \hspace{1.5pt} y_{d\minus1}, \;  \omega \in \Omega \\
					& \;\overline{\vartheta}^{\omega n}_{d\plus1} \geq \sum\limits_{s \in \mathcal{S}} \mu^{\omega n}_{s,d\plus1} \overline{\theta}^{n}_{s,d\plus1} + \mathrm{M_y} ||\gamma^{\omega n}_{d\plus1}|| + \mathrm{M_b} ||\zeta^{\omega n}_{d\plus1}||, \;  \omega \in \Omega,  n \in N \\
					& \; \sum\nolimits_{s \in \mathcal{S}} \mu^{\omega n}_{s,d\plus1} y^{n}_{s,d} = y^\omega_d + \gamma^{\omega n}_{d\plus1}, \;  \omega \in \Omega,  n \in \mathcal{M} \\ 
					& \; \sum\nolimits_{s \in \mathcal{S}} \mu^{\omega n}_{s,d\plus1} b^{n}_{s,d\plus1} = b^{mn}_{d\plus1} + \zeta^{\omega m}_{d\plus1}, \;  \omega \in \Omega,  n \in \mathcal{M} \\ 
					& \; \sum\nolimits_{s \in \mathcal{S}} \mu^{\omega n}_{s,d\plus1} = 1, \;  \omega \in \Omega,  n \in \mathcal{M}. 
		\end{aligned}
	\end{equation}
	\begin{lemma}\label{lemma:UPbf_Lemma}
		Solving problem $\mathbf{UP}^{m}_{d}(y_{d-1},b^{lm}_d)$ yields an optimal solution $\overline{y}_d$ and an optimal objective $\overline{\theta}_d$ which is a valid upper bound on the objective of $\mathbf{SP}^{m}_{d}(x,y_{d-1},b^{lm}_d)$.
	\end{lemma}
	\begin{proof}
		This lemma can be proved with the same argument of Lemma~\ref{lemma:LP_Lemma}. \qed
	\end{proof}

	\subsection{Enhanced SDDP algorithm}
	The extended SDDP is presented in Algorithm~\ref{alg:SDDPv1}. In the following, we show the convergence proof of the algorithm.

 \begin{algorithm}[!htb]
	\caption{\texttt{Enhanced SDDP algorithm}}
	\label{alg:SDDPv1}
	choose tolerance $\delta > 0$ and set $k:=0$\;
	set $\underline{\theta}^k:=-\infty$, $\overline{\theta}^k:=\infty$, and set $\underline{\theta}^{m}_1:=-\infty$, $\overline{\theta}^{m}_1:=\infty$, ${\scriptstyle m \in \mathcal{M}}$\; 
	\Repeat{$\overline{\theta}^{k}-\underline{\theta}^{k}\leq\delta$}{
		set $k:=k+1$, $d:=1$, and $\hat{d}:=\mathbf{d}$\;
		\tcc{Forward Pass}
		set $m_1:=$\raisebox{.75pt}{\scalebox{.90}{$\mathrm{arg} \,\underset{m \in \mathcal{M}}{\mathrm{max}}  \left\{ \left( \overline{\theta}^m_1 \hspace{-3pt} - \underline{\theta}^m_1 \hspace{-1pt} \right) \,\middle|\, \pi^m_1 > 0 \right\}$}}\;
		solve $\mathbf{LP}^{m_1}_{1}$, get $\underline{y}_1$ and $\{\underline{\vartheta}^{n}_{2}, {\scriptstyle n \in \mathcal{M}}\}$\;
		solve $\mathbf{UP}^{m_1}_{1}(\underline{y}_1)$, get $\{\overline{\underline{\vartheta}}^{n}_{2}, {\scriptstyle n \in \mathcal{M}}\}$\;
		\Repeat{$\hat{d}=d-1$ $\mathbf{or}$ $d=\mathbf{d}$}{
			set $d:=d+1$\;
			set $(\omega_{d\minus1},m_d):=$\raisebox{.75pt}{\scalebox{.90}{$\mathrm{arg} \underset{\omega \in \Omega, n \in \mathcal{M}}{\mathrm{max}} \left\{ \left( \overline{\underline{\vartheta}}^{\omega n}_d \hspace{-3pt} - \underline{\vartheta}^{\omega n}_d \right) \;\middle|\; \pi^{m_{d\minus1} n}_d > 0 \right\}$}}\;
			\uIf{\raisebox{.75pt}{\scalebox{.90}{$\left( \overline{\underline{\vartheta}}^{\omega_{d\minus1} m_d}_d \hspace{-3pt} - \underline{\vartheta}^{\omega_{d\minus1} m_d}_d \hspace{-1pt} \right)$}}$ \, > \tfrac{\delta (\mathbf{d}-d+1)}{\mathbf{d}-1} $}{
				solve $\mathbf{LP}^{m_d}_{d}(x,\underline{y}^{\omega_{d\minus1}}_{d\minus1},b^{m_{d\minus1}m_d}_d)$, get $\underline{y}_d$ and $\{\underline{\vartheta}^{\omega n}_{d\plus1}, {\scriptstyle \omega \in \Omega, n \in \mathcal{M}}\}$\;
				solve $\mathbf{UP}^{m_d}_{d}(\underline{y}^{\omega_{d\minus1}}_{d\minus1},b^{m_{d\minus1}m_d}_d,\underline{y}_{d})$, get $\{\overline{\underline{\vartheta}}^{\omega n}_{d\plus1}, {\scriptstyle \omega \in \Omega, n \in \mathcal{M}}\}$\;
				}\Else{
				$\hat{d}:=d-1$\;	
				}
		}
		\tcc{Backward Pass}
		\For{$d=\hat{d},\hat{d}\raisebox{1.5pt}{\scalebox{.7}{$-$}}1,..,2$}{
			\For{$m \in \mathcal{M}$}{
				solve $\mathrm{LP}^{m}_{d}(x,\underline{y}^{\omega_{d\minus1}}_{d\minus1},b^{m_{d\minus1}m}_d)$, get $\underline{\theta}^{m}_{d}$, $\underline{\lambda}^{m}_{d}$, $\underline{\sigma}^{m}_{d}$, $\underline{\nu}^{m}_{d}$\;
				solve $\mathrm{UP}^{m}_{d}(\underline{y}^{\omega_{d\minus1}}_{d\minus1},b^{m_{d\minus1}m}_d)$, get $\overline{\theta}^{m}_{d}$\;
				\For{$l \in \mathcal{M}$}{
					add cuts $\left(\underline{\theta}^{m}_{d},\underline{\lambda}^{m}_{d},\underline{\sigma}^{m}_{d},\underline{\nu}^{m}_{d},x,\underline{y}^{\omega_{d\minus1}}_{d\minus1},b^{m_{d\minus1}m}_d\right)$ to $\mathbf{LP}^{l}_{d\minus1}$\;
					add cuts $\left(\overline{\theta}^{m}_{d},\underline{y}^{\omega_{d\minus1}}_{d\minus1},b^{m_{d\minus1}m}_d\right)$ to $\mathbf{UP}^{l}_{d\minus1}$\;
				}
			}
		}
		\For{$m \in \mathcal{M}$}{
			solve $\mathrm{LP}^{m}_{1}(x)$ and $\mathrm{UP}^{m}_{1}$ to get $\underline{\theta}^{m}_{1}$, $\underline{\lambda}^{m}_{1}$, and $\overline{\theta}^{m}_{1}$\;
		}
		set $\underline{\theta}^{k} := \sum\limits_{m \in \mathcal{M}}\pi^m_1\underline{\theta}^{m}_{1}$, $\underline{\lambda}^{k} := \sum\limits_{m \in \mathcal{M}}\pi^m_1\underline{\lambda}^{m}_{1}$, and $\overline{\theta}^{k}  := \sum\limits_{m \in \mathcal{M}}\pi^m_1\overline{\theta}^{m}_{1}$ \;
	}
\end{algorithm}

	\subsubsection*{Convergence proof of the extended SDDP algorithm}

	\begin{theorem}\label{theorem:conv_sddp2}
		For given convergence tolerance $\delta>0$, Algorithm~\ref{alg:SDDPv1} converges to an $\delta$-optimal solution in a finite number of iterations.
	\end{theorem}
	\begin{proof} \hspace{1pt} \\
		Notice that the forward pass of Algorithm~\ref{alg:SDDPv1} is equivalent to the forward pass of Algorithm~\ref{alg:SDDPv0}. In the backward pass, Algorithm~\ref{alg:SDDPv0} builds cuts only on the path $\{m_d, d=1,..,\hat{d}\}$ chosen during the foward pass. The backward pass of Algorithm~\ref{alg:SDDPv1}, instead, solves each problem $\mathbf{LP}^m_d$ and $\mathbf{UP}^m_d$ for each $m \in \mathcal{M}$ and then add each generated cut to each problem $\mathbf{LP}^l_{d\minus1}$ and $\mathbf{UP}^l_{d\minus1}$ for each $l \in \mathcal{M}$. Algorithm~\ref{alg:SDDPv1} adds at least all the cuts that would be added by Algorithm~\ref{alg:SDDPv0} in the backward pass of the same iteration. Adding more valid cutting planes than Algorithm~\ref{alg:SDDPv0} does not affect the convergence properties of the SDDP algorithm. It follows that also Algorithm~\ref{alg:SDDPv1} converges to an $\delta$-optimal solution in a finite number of iterations 
		\qed
	\end{proof}

\section{Adaptive Benders decomposition}\label{sec:Benders}
We use the relaxed master problem $\mathbf{RMP}$ defined in~\eqref{eq:RMP} in Adaptive Benders decomposition algorithm to iteratively solve the full problem $\mathbf{MP}$ defined in~\eqref{eq:opt_init} up to an $\epsilon$-optimal solution. 
At each iteration $w$ we solve problem $\mathbf{RMP}$ and obtain a set of decisions $\mathbf{x}^{w}$. For given master decisions $\mathbf{x}^{w}$ we solve the set $I$ of subproblems $\mathbf{SP}(x^w_i)$ up to an $\delta$-optimal solution and obtain a valid lower bound $\underline{\theta}^w_i$, a valid subgradient $\underline{\lambda}^w_i$ w.r.t. $x$, and a valid upper bound $\overline{\theta}^w_i$. Then, a set of cutting planes are added to the $\mathbf{RMP}$ at points $\mathbf{x}^{w}$. At each iteration $w$, the Benders algorithm computes a valid lower bound $\mathrm{L}^w$ and a valid upper bound $\mathrm{U}^w$ on the optimal objective of problem $\mathbf{MP}$. The algorithm stops when $\mathrm{U}^w-\mathrm{L}^w \leq \epsilon$. The key of Adaptive Benders decomposition is to exploit the subproblem structure and conduct cut sharing. The adaptive oracles were introduced for problems where the following conditions hold:

\begin{cond}
   $\mathbf{SP}(x^w_i)$ is convex w.r.t. the vector $x^{w}_i$, and $\mathbf{SP}(x^w_i)$ is a decreasing function of the elements of $x^{w}_i$.  \label{propty: g(x,q)}
\end{cond}

The convexity is immediate consequence of $\mathbf{SP}(x^w_i)$ being a minimisation linear program and the monotonicity properties hold if, for example, $A_d, B_d$ and $y_d$ are non-negative. 

Once one or more subproblems have been solved at a collection of points, this information can be used by the adaptive oracles to generate valid bounds at different points for all subproblems. We refer to \cite{Mazzi2020} for the mathematical definition and the proof of the properties of the adaptive oracles.

The adaptive oracles provide bounds for a subproblem at a new solution point without having to solve it exactly, and this reduces the computational cost compared to standard Benders decomposition. The process is shown in Algorithm \ref{alg:Bend}. In iteration \(w\), when subproblem \(\hat{i}\) is solved at the point \(\hat{x}^w_{\hat{i}}\) using Algorithm \ref{alg:SDDPv1}, the algorithm returns the optimal value \(\theta_{\hat{i}}^{w} = g(\hat{x}_{\hat{i}}^{w})\), and the subgradients \(\lambda_{\hat{i}}^{w}\) with respect to \(x_{\hat{i}}^{w}\). Then the solution vector \((\hat{x}_{\hat{i}}^{w}, \theta_{\hat{i}}^{w}, \lambda_{\hat{i}}^{w})\) is added to the collection \(\mathcal{Z}\) of solution vectors. Then, using the information in \(\mathcal{Z}\), the adaptive oracles generate valid bounds for all subproblems: the oracles are called for each subproblem \(i\) at the current solution point \(\hat{x}^w_{i}\) and return the values \(\underline{\theta}_{i}^{w}, \overline{\theta}_{i}^{w}\), and \(\lambda_{i}^{w}\) with the properties:
\begin{proper}$\underline{\theta}_{i}^{w} + \lambda^{w\top}_{i}(x^w_i - \hat{x}^w_{i}) \le g(x^w_i),~~ \forall x^w_i~~~ \text{and} ~~~ g(\hat{x}^w_{i}) \le \overline{\theta}_{i}.
$ \label{propty: cuts}
\end{proper}

The RMP in adaptive Benders is the same as in standard Benders, except that the exact cuts of standard Benders are replaced by the approximate cuts in, which use the quantities supplied by the adaptive oracles. 

\subsection{Adaptive Benders decomposition algorithm}
The Adaptive Benders decomposition is presented in Algorithm~\ref{alg:Bend}.
\begin{algorithm}[!htb]
    \caption{\texttt{Adaptive Benders decomposition}}
    \label{alg:Bend}
    choose $\delta>0$, $\epsilon > 2 \sum_{i \in \mathcal{I}} \pi_i \hspace{1pt} \delta$, and set $w:=0$, $\mathrm{L}^w:=-\infty$ and $\mathrm{U}^w:=\infty$\;
    solve subproblem exactly at $\underline{x}$ and obtain $\theta$, and $\lambda$; set $\mathcal{S}:=\{(\underline{x},\theta,\lambda)\}$\;
    \Repeat{$\mathrm{U}^{w}-\mathrm{L}^{w} \leq \epsilon$}{
        set $w:=w+1$\;
        solve RMP and obtain $\beta_{i}^{w}$ and $\mathbf{x}^{w}$\;
        set $\mathrm{L}^{w}:=\sum_{i \in \mathcal{I}}\pi_{i}(c^{\top}x_{i}^{w} +\beta_{i}^{w})$\;
        \For{$i \in \mathcal{I}$}{
            call adaptive oracles at $x_{i}^{w}$ and obtain $\underline{\theta}_{i}^{w}$, $\overline{\theta}_{i}^{w}$, and $\underline{\lambda}_{i}^{w}$\;
        }
        set $\xi:=0$\;  
        \Repeat{$\xi=1$ \textbf{or} }{
            \If{$\max \pi_i(\overline{\theta}_{i}^{w}-\underline{\theta}_{i}^{w})=0$}{\textbf{break}}
            set $\hat{i}:= \arg\max_{i \in \mathcal{I}} \pi_i (\overline{\theta}_{i}^{w} - \underline{\theta}_{i}^{w})$\;
            solve subproblem exactly at $x_{\hat{i}}^{w}$ using Algorithm \ref{alg:SDDPv1} and obtain $\theta_{\hat{i}}^{w}$, $\lambda_{\hat{i}}^{w}$\;
            \If{$\theta_{\hat{i}}^{w}>\underline{\theta}_{\hat{i}}^{w}$}{
                set $\xi:=1$\;
                update $\mathcal{Z}:=\mathcal{Z}\cup\{(x_{\hat{i}}^{w},\theta_{\hat{i}}^{w},\lambda_{\hat{i}}^{w})\}$\;
            }
        }
            \For{$i \in \mathcal{I}$}{
                call adaptive oracles at $\hat{x}_{i}^{w}$ and update $\underline{\theta}_{i}^{w}$, $\overline{\theta}_{i}^{w}$, $\underline{\lambda}_{i}^{w}$\;
            }
        set $\mathrm{U}^{w}:=\min(\mathrm{U}^{w-1}, \sum_{i \in \mathcal{I}}\pi_{i}(c^{\top}x_{i}^{w} +\overline{\theta}_{i}^{w}))$\;
    }
\end{algorithm}

\subsubsection*{Convergence proof of the Adaptive Benders decomposition}
The convergence proof relies on $\mathcal{X}$ being a compact set and on $g(x)$ defined in~\eqref{eq:multistage_g} having bounded subgradients. Therefore, we assume that it exists a finite $\alpha_x$ such that $||\underline{\lambda}(x)||\leq \alpha_x$ for all $x \in \mathcal{X}$, where $\underline{\lambda}(x)$ is a valid subgradient w.r.t. $x$ of such that 
$$ \underline{\theta}(x) + \underline{\lambda}(x)^\top \left(\hat{x} - x \right),  \hat{x} \in \mathcal{X} \quad \text{and} \quad \underline{\theta}(x)\leq g(x) \leq \underline{\theta}(x) + \delta, \quad  x \in \mathcal{X}. $$

\begin{lemma}\label{lemma:bend_0}
	Let $x^w_i$ and $\beta^w_i$ be the optimal solution of $\mathbf{RMP}$ at iteration $w$ for subproblem $i$. Then, let $\overline{\theta}^w_i$ and $\underline{\theta}^w_i$ be valid upper and lower bound on the true objective $g(x^w_i)$ such that $\overline{\theta}^w_i-\underline{\theta}^w_i \leq \delta$. If it exists $\hat{w}\in\{1,..,w-1\}$ such that 
	$$||x^w_i-x^{\hat{w}}_i||\leq \tfrac{\epsilon-2\delta\sum_{i \in \mathcal{I}}  \pi_i }{2 \alpha_x \sum_{i \in \mathcal{I}} \pi_i }, \; \text{where} \; \epsilon > 2\delta \textstyle\sum\limits_{i \in \mathcal{I}} \pi_i $$
	it follows that
	$\overline{\theta}^w_i-\beta^w_i \leq \tfrac{\epsilon}{\sum_{i \in \mathcal{I}} \pi_i}$.	
\end{lemma}
\begin{proof} \hspace{1pt} \\
	Given that $g(x)$ has bounded subgradients it follows that $||\underline{\lambda}^{\hat{w}}_i|| \leq \alpha_x$ and
	\begin{equation}\label{eq:lemma_bend_1}
		\beta^w_i \geq \underline{\theta}^{\hat{w}}_i - \alpha_x ||x^w_i-x^{\hat{w}}_i||,
	\end{equation}	
	since $\beta^w_i \geq \underline{\theta}^{\hat{w}}_i + \underline{\lambda}^{{\hat{w}}\top}_i \hspace{-2pt} \left(x^w_i-x^{\hat{w}}_i\right)$. Then, the definition of $\underline{\theta}^w_i$ computed by Algorithm~\ref{alg:Bend} solving $\mathbf{SP}(x^w_i)$ up to a $\delta$-optimal solution leads to
	\begin{equation}\label{eq:lemma_bend_2}
		\begin{aligned}
			\overline{\theta}^w_i & \leq g(x^w_i) + \delta\\
			                      & \leq g(x^{\hat{w}}_i) + \alpha_x||x^w_i-x^{\hat{w}}_i|| + \delta\\
			                      & \leq \overline{\theta}^{\hat{w}}_i + \alpha_x||x^w_i-x^{\hat{w}}_i|| + \delta.\\
		\end{aligned}
	\end{equation}
	The first inequality holds since $g(x^w_i)+\delta$ is a valid upper bound on $\overline{\theta}^w_i$,
	the second inequality holds since $g(x)$ has bounded subgradients, and the third inequality holds since $\overline{\theta}^{\hat{w}}_i$ is a valid upper bound on $g(x^{\hat{w}}_i)$.
	Combining equation~\eqref{eq:lemma_bend_1} with equation~\eqref{eq:lemma_bend_2} leads to
	\begin{equation}\label{eq:lemma_bend_3}
		\begin{aligned}
			\overline{\theta}^w_i - \beta^w_i & \leq \Big( \overline{\theta}^{\hat{w}}_i + \alpha_x||x^w_i-x^{\hat{w}}_i|| + \delta \Big) - \Big( \underline{\theta}^{\hat{w}}_i - \alpha_x ||x^w_i-x^{\hat{w}}_i|| \Big)  \\
			& \leq 2 \delta + 2 \alpha_x \tfrac{\epsilon-2\delta\sum_{i \in \mathcal{I}} \pi_i }{2 \alpha_x \sum_{i \in \mathcal{I}}\pi_i }   \\
			& \leq \tfrac{\epsilon}{\sum_{i \in \mathcal{I}} \pi_i}.  \\
		\end{aligned}
	\end{equation} \qed
\end{proof}

\begin{lemma}\label{lemma:bend_1}
	There exists a finite number of cutting planes that can be added to $\mathbf{RMP}$, one at a time, for subproblem $i$, such that 
	$$||x^w_i-x^{\hat{w}}_i|| > \tfrac{\epsilon-2\delta\sum_{i \in \mathcal{I}} \pi_i }{2 \alpha_x \sum_{i \in \mathcal{I}}\pi_i }, \quad  {\hat{w}}=1,..,w-1$$
\end{lemma}
\begin{proof} \hspace{1pt} \\
	Observe that $\tfrac{\epsilon-2\delta\sum_{i \in \mathcal{I}} \pi_i }{2 \alpha_x \sum_{i \in \mathcal{I}}\pi_i }$ is a positive and finite number since $\epsilon > 2\delta \textstyle\sum_{i \in I} \pi_i$ by assumption, and $\alpha_x$ is a finite constant. Given that $\mathcal{X}$ is a compact set, it follows that there exists a finite number of cutting planes that can be added to $\mathbf{RMP}$ for subproblem $i$ such that $||x^w_i-x^{\hat{w}}_i|| > \tfrac{\epsilon-2\delta\sum_{i \in \mathcal{I}} \pi_i }{2 \alpha_x \sum_{i \in \mathcal{I}}\hspace{-2pt}\pi_i }$ for each ${\hat{w}}=1,..,w-1$.
	\qed
\end{proof}

\begin{theorem}\label{theorem:conv_bend}
	For given convergence tolerance $\epsilon>2\delta \textstyle\sum_{i \in \mathcal{I}}\pi_i$, Algorithm~\ref{alg:Bend} converges to an $\epsilon$-optimal solution in a finite number of iterations.
\end{theorem}
\begin{proof} \hspace{1pt} \\
	By Lemma~\ref{lemma:bend_1}, we know that there exists a finite number of cutting planes that can be added to $\mathbf{RMP}$, one at the time, for subproblem $i$, such that $||x^{w}_i-x^{\hat{w}}_i|| >\tfrac{\epsilon-2\delta\sum_{i \in \mathcal{I}} \hspace{-2pt}\pi_i }{2 \alpha_x \sum_{ \in \mathcal{I}}\hspace{-2pt}\pi_i }$ for each ${\hat{w}}=1,..,w-1$. Given that the number $|\mathcal{I}|$ of subproblems is finite, it follows that after a finite number of iterations Algorithm~\ref{alg:Bend} reaches a solution $\mathbf{x}^{w}$ such that
	\begin{equation*}
		 ||x^w_i-x^{{\hat{w}}_i}_i||\leq\tfrac{\epsilon-2\delta\sum_{i \in \mathcal{I}} \hspace{-2pt}\pi_i }{2 \alpha_x \sum_{i \in \mathcal{I}}\hspace{-2pt}\pi_i }, \quad {\hat{w}}_i \in \{1,..,w-1\}, \quad  i \in \mathcal{I}
	\end{equation*}
	and hence $\overline{\theta}^w_i - \beta^w_i \leq \tfrac{\epsilon}{\sum_{i \in \mathcal{I}}\pi_i},$ $ i \in \mathcal{I}$. It follows that
	
	\begin{equation}\label{eq:therem_bend}
		\begin{aligned}
			\mathrm{U}^{w} - \mathrm{L}^{w} & \leq \left( f(\mathbf{x}^{w}) + \sum_{i \in \mathcal{I}} \hspace{-1pt} \pi_i \overline{\theta}^w_i \right) - \left( f(\mathbf{x}^{w}) + \sum_{i \in \mathcal{I}} \hspace{-1pt} \pi_i \beta^w_i \right) \\
			                      & = \sum_{i \in \mathcal{I}} \hspace{-1pt} \pi_i \left( \overline{\theta}^w_i - \beta^w_i  \right)\\
			                      & \leq \epsilon.\\
		\end{aligned}
	\end{equation}\qed
\end{proof}

\section{Case Study}\label{sec:casestudy}
We test the proposed algorithms on a power system stochastic planning problem with short-term (wind power) and long-term (energy demand level) uncertain parameters. The algorithms are implemented in \textsc{Julia} 1.11.6. A MacBook Pro with an Apple M4 Max chip and 64 GB of RAM is used for running the code. The optimisation models are implemented in \textsc{JuMP} \cite{jump} and solved with \textsc{Gurobi} 12.0.2.

\subsection{Investment planning model}

We consider a power system investment planning problem with a time horizon of 15 years. The investment problem has 7 decision nodes: one refers to decisions to be taken in the first stage, two to decisions in 5 years time, and four to decisions in 10 years time. At each node, we also compute the cost of operating the system for the following 5 years for a given installed capacity. We consider a construction time of 5 years, so new assets installed in the first stage will only be available in 5 and 10 years, and new capacity installed in 5 years will only be available in 10 years. We model a set of $\mathcal{P}$ technologies: 6 thermal units, 1 storage unit, and 2 renewable generation units.

We formulate the stochastic investment planning problem as~\eqref{eq:opt_init}
where $\mathcal{I}$ is the set of stochastic decision nodes, each associated with a probability $\pi_i$. The vector of coefficients $x_i$ is given by
\begin{equation*}
	x_i = \left(\left\{x^{acc}_{pi},  p \in \mathcal{P} \right\},\nu^{D}_i\right), \quad  i \in \mathcal{I}, 
\end{equation*}
where $x^{acc}_{pi}$ is the accumulated capacity of technology $p$ at node $i$. Parameter $\nu^{D}_i$ is the relative level of energy demand. 

The cost for operating the system for a given vector $x_i$ is obtained by solving a set of multistage stochastic and linear programs like~\eqref{eq:sp_init_d} via Algorithm~\ref{alg:SDDPv1}. The uncertain parameter $b^{lm\omega}_d$ represents the 24-hour trajectory of wind power production (in per unit) during day $d$ that ends at state $m$ via scenario $\omega$ given that day $d-1$ ended at state $l$. We consider 4 slices, one per each season, each of which has $\mathbf{d}=7$ stages, $\mathbf{m}=5$ states, and $\bm{\omega}=3$ scenarios. 

At each stage $d$, we solve a 24-hour economic dispatch for given wind power productions $\{b^{lm\omega}_d, \omega=1,..,\bm{\omega} \}$ and for a given value of decision $y_{d\minus1}$ fixed during $d-1$. Only a subset of $y_{d\minus1}$ actually influences the problem at stage $d$, and we define such subset as $\hat{y}_{d\minus1}=C_d \, y_{d\minus1}$. In our problem, $\hat{y}_{d\minus1}$ is a vector that includes the level of generation of slow thermal units (coal, coal\&CCS, and nuclear) and the level of the energy left in a storage unit at the end of day $d-1$. The first 3 elements of $\hat{y}_{d\minus1}$ enforce ramping limitations during the first hour of day $d$, while the last element imposes the energy conservation to the storage unit. 

\begin{figure}[htbp!] \label{fig:DAM_Scenarios}  
  \centering
      \includegraphics[width=1.00\columnwidth]{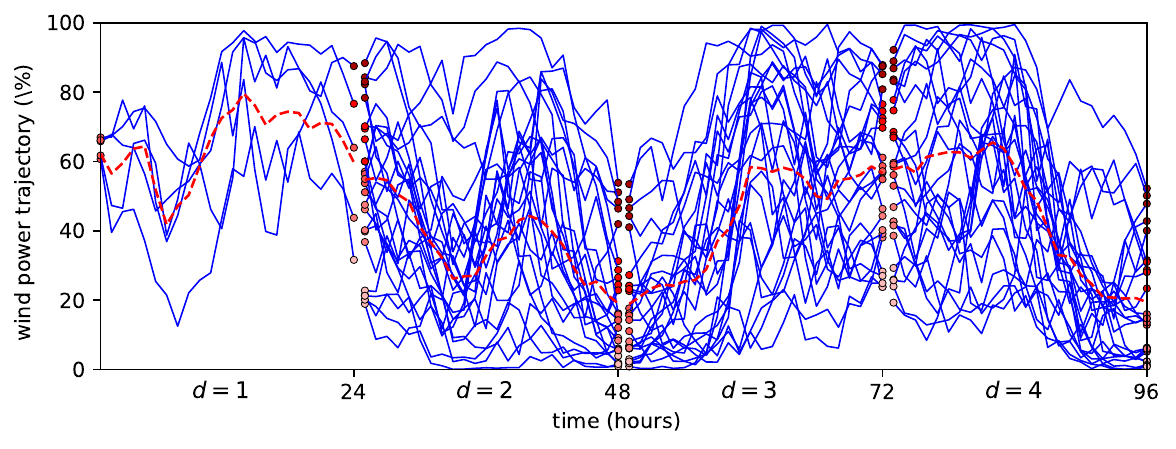}
	  \caption{Example of short-term wind power uncertainty for 3 days. The grey continuous lines represent the wind trajectories $b^{lm\omega}_d$, and the squares at the end of each day indicate the state $m$ they represent. Some of the trajectories are highlighted to show some of the possible paths.}
\end{figure}

As a benchmark, we solve the investment planning problem~\eqref{eq:opt_init} ignoring the short-term uncertainty in wind power, i.e., the subproblem is a standard economic dispatch with $\mathbf{d}=7$ stages.
To make sure the deterministic subproblem is still able to capture the variability of wind power generation, we sample 500 possible wind power trajectories, and we select the $\bm{\iota}$ more representative ones via a \texttt{kmeans} clustering algorithm. Each trajectory $\{ b^{\iota}_d, d=1,..,\bm{d} \}$ is associated with a probability $\pi^{\iota}$ such that $\sum_{\iota=1}^{\bm{\iota}} \pi^\iota = 1$. To obtain the cost for operating the system for a given vector $x_i$ with the deterministic subproblem, we therefore solve $\bm{\iota}$ economic dispatch and we weight each of them with the associated probability $\pi^{\iota}$. 



\subsection{Results}\label{sec:results}

We solve the stochastic investment problem with short-term uncertainty in the operational subproblem, and we compare it with its deterministic alternative, but still with long-term uncertainty.

Table \ref{table:optimal_objectives} shows the difference in the optimal objectives between the two cases. The deterministic case yields a 10.7\% lower optimal objective compared with the stochastic version. This is because, in the deterministic case, the model has perfect weekly foresight of the wind capacity factor, which is not the case in the real world. This leads the model to obtain a lower (and underestimated) total cost than its stochastic counterpart.

\begin{table}[htb!]
	\centering
	\caption{Optimal objectives ($10^9\pounds$).}
		\begin{tabular}{ccc}
		\toprule
		\multirow{2}{*}{version} & lower & upper \\
								& bound & bound \\ \midrule
		deterministic            & 119.3 & 119.7 \\
		stochastic               & 133.6 & 134.3 \\ \bottomrule
	\end{tabular}
    \label{table:optimal_objectives}
\end{table}

\begin{table}[htb!]
	\centering
	\caption{Optimal investments (GW) at nodes 1,2, and 3 when the subproblem is deterministic ($\bm{\iota}=100$).}
	\label{tab:det_inv}
	\begin{tabular}{cccccc}
	\toprule
			  &           & \multirow{3}{*}{\begin{tabular}{@{}c@{}}historical\\capacity\end{tabular}} & \multicolumn{3}{c}{newly installed capacity} \\
	tech. $p$ & type      &            & node 1        & node 2        & node 3       \\
			  &           &            & (present)     & (5 years)     & (5 years)    \\ \midrule
	coal      & thermal   & 8.8       & 0.0           & 0.0           & 0.0          \\
	coal\&CCS & thermal   & 1.6        & 0.0           & 0.0           & 0.0          \\
	OCGT      & thermal   & 4.0        & 0.0           & 2.0           & 0.0          \\
	CCGT      & thermal   & 13.6       & 1.0           & 3.4           & 5.2          \\
	diesel    & thermal   & 0.8        & 0.5           & 0.0           & 4.5          \\
	nuclear   & thermal   & 4.0        & 16.8          & 1.0           & 5.1          \\
	lithium   & storage   & 0.4       & 75.9          & 0.0           & 25.2         \\
	wind      & renewable & 4.8        & 53.8          & 2.2           & 14.2         \\
	solar     & renewable & 4.4        & 0.0           & 0.0           & 0.0          \\ \bottomrule
	\end{tabular}
    \label{table: deterministic_solution}
\end{table}

\begin{table}[htb!]
	\centering
	\caption{Optimal investments (GW) at nodes 1,2, and 3  when the subproblem is stochastic.}
	\label{tab:sto_inv}
	\begin{tabular}{cccccc}
	\toprule
			  &           & \multirow{3}{*}{\begin{tabular}{@{}c@{}}historical\\capacity\end{tabular}} & \multicolumn{3}{c}{newly installed capacity} \\
	tech. $p$ & type      &            & node 1        & node 2        & node 3       \\
			  &           &            & (present)     & (5 years)     & (5 years)    \\ \midrule
	coal      & thermal   & 8.8       & 0.0           & 0.0           & 0.0          \\
	coal\&CCS & thermal   & 1.6        & 0.0           & 0.0           & 0.0          \\
	OCGT      & thermal   & 4.0        & 0.0           & 2.8           & 0.4          \\
	CCGT      & thermal   & 13.6       & 2.0           & 2.4           & 4.7          \\
	diesel    & thermal   & 0.8        & 2.2           & 0.0           & 4.2          \\
	nuclear   & thermal   & 4.0        & 18.7          & 1.5           & 7.8          \\
	lithium   & storage   & 0.4       & 72.2          & 14.0          & 0.0          \\
	wind      & renewable & 4.8        & 39.6          & 1.5           & 4.8          \\
	solar     & renewable & 4.4        & 0.0           & 0.0           & 0.0          \\ \bottomrule
	\end{tabular}
    \label{table:stochastic_solution}
\end{table}

Tables \ref{table: deterministic_solution} and \ref{table:stochastic_solution} present the investment decisions in the first three nodes in the two cases, respectively. By comparing the optimal investments, we notice that the overall investment mix is similar in both cases. The stochastic case has higher investment in nuclear, which is used to balance the wind generation volatility. Interestingly, we find that the total investment in lithium batteries is lower in the stochastic case. This is because the deterministic model makes more investment in wind in node 1 and less investment in nuclear than the stochastic case. Hence, the deterministic case requires more storage to balance the wind uncertainty.

\section{Conclusions}\label{sec:conclusions}
This paper proposes the first algorithm to solve multistage stochastic programmes with block-separable multistage recourse. An example of such a problem is MHSP with long-term and short-term uncertainty both revealed at multiple stages. The proposed algorithm has two parts: (1) Adaptive Benders decomposition to decompose the whole problem into a reduced master problem and independent blocks of subproblems, and (2) an enhanced SDDP to solve each independent subproblem with multistage uncertainty. The algorithm is applied to solve a power system planning problem with long-term and short-term uncertainty. The case study results show that (1) the proposed algorithm can efficiently solve this type of problem, (2) deterministic wind modelling underestimates the objective function, and (3) stochastic modelling of wind leads to different investment decisions. Future research includes stabilisation of such algorithms and further cut sharing. Note that in this paper, we only consider short-term and long-term uncertainty in the right-hand-side parameters, but it is straightforward to extend the approach to uncertainty in cost coefficients.

\section*{CRediT author statement}
\textbf{Nicolò Mazzi:} Conceptualisation, Methodology, Data curation, Software, Validation, Investigation, Formal analysis, Writing - original draft. \textbf{Ken McKinnon:} Conceptualisation, Methodology, Supervision, Writing - original draft. \textbf{Hongyu Zhang:} Conceptualisation, Methodology, Software, Visualisation, Writing - original draft.

\bibliographystyle{spmpsci}  
\bibliography{benders_ao_sddp_v0}      

\begin{thebibliography}{10}
\providecommand{\url}[1]{{#1}}
\providecommand{\urlprefix}{URL }
\expandafter\ifx\csname urlstyle\endcsname\relax
  \providecommand{\doi}[1]{DOI~\discretionary{}{}{}#1}\else
  \providecommand{\doi}{DOI~\discretionary{}{}{}\begingroup \urlstyle{rm}\Url}\fi

\bibitem{Backe2022}
Backe, S., Skar, C., del Granado, P.C., Turgut, O., Tomasgard, A.: {EMPIRE: an open-source model based on multi-horizon programming for energy transition analyses}.
\newblock SoftwareX \textbf{17}, 100877 (2022).
\newblock \doi{10.1016/j.softx.2021.100877}

\bibitem{Benders1962}
Benders, J.F.: {Partitioning procedures for solving mixed-variables programming problems}.
\newblock Numerische Mathematik \textbf{4}(1), 238--252 (1962).
\newblock \doi{10.1007/BF01386316}

\bibitem{Birge1985DecompositionPrograms}
Birge, J.R.: {Decomposition and Partitioning Methods for Multistage Stochastic Linear Programs}.
\newblock Operations Research \textbf{33}(5), 989--1007 (1985).
\newblock \doi{10.1287/OPRE.33.5.989}

\bibitem{Downward2020StochasticUncertainty}
Downward, A., Dowson, O., Baucke, R.: {Stochastic dual dynamic programming with stagewise-dependent objective uncertainty}.
\newblock Operations Research Letters \textbf{48}(1), 33--39 (2020).
\newblock \doi{10.1016/J.ORL.2019.11.002}

\bibitem{jump}
Dunning, I., Huchette, J., Lubin, M.: {JuMP: A Modeling Language for Mathematical Optimization}.
\newblock SIAM Review \textbf{59}(2), 295--320 (2017).
\newblock \doi{10.1137/15M1020575}

\bibitem{Durakovic2024DecarbonizingSectors}
Durakovic, G., Zhang, H., Knudsen, B.R., Tomasgard, A., del Granado, P.C.: {Decarbonizing the European energy system in the absence of Russian gas: Hydrogen uptake and carbon capture developments in the power, heat and industry sectors}.
\newblock Journal of Cleaner Production \textbf{435}, 140473 (2024).
\newblock \doi{10.1016/J.JCLEPRO.2023.140473}

\bibitem{Fullner2021StochasticVariants}
F{\"{u}}llner, C., Rebennack, S.: {Stochastic dual dynamic programming and its variants}.
\newblock Available at Optimization Online  (2021).
\newblock \urlprefix\url{http://www.optimization-online.org/DB_FILE/2021/01/8217.pdf}

\bibitem{Homem-De-Mello2011SamplingScheduling}
Homem-De-Mello, T., De~Matos, V.L., Finardi, E.C.: {Sampling strategies and stopping criteria for stochastic dual dynamic programming: A case study in long-term hydrothermal scheduling}.
\newblock Energy Systems \textbf{2}(1), 1--31 (2011).
\newblock \doi{10.1007/S12667-011-0024-Y/METRICS}

\bibitem{Kaut2014}
Kaut, M., Midthun, K.T., Werner, A.S., Tomasgard, A., Hellemo, L., Fodstad, M.: {Multi-horizon stochastic programming}.
\newblock Computational Management Science \textbf{11}(1-2), 179--193 (2014).
\newblock \doi{10.1007/s10287-013-0182-6}

\bibitem{King2024ModelingEdition}
King, A.J., Wallace, S.W.: {Modeling with Stochastic Programming: Second Edition}.
\newblock Springer Series in Operations Research and Financial Engineering \textbf{Part F2932}, 1--197 (2024).
\newblock \doi{10.1007/978-3-031-54550-4/COVER}

\bibitem{Lara2020}
Lara, C.L., Siirola, J.D., Grossmann, I.E.: {Electric power infrastructure planning under uncertainty: stochastic dual dynamic integer programming (SDDiP) and parallelization scheme}.
\newblock Optimization and Engineering \textbf{21}(4), 1243--1281 (2020).
\newblock \doi{10.1007/s11081-019-09471-0}

\bibitem{Leclere2020ExactDuality}
Lecl{\`{e}}re, V., Carpentier, P., Chancelier, J.P., Lenoir, A., Pacaud, F.: {Exact Converging Bounds for Stochastic Dual Dynamic Programming via Fenchel Duality}.
\newblock SIAM Journal on Optimization \textbf{30}(2), 1223--1250 (2020).
\newblock \doi{10.1137/19M1258876}

\bibitem{Louveaux1986MultistageRecourse}
Louveaux, F.V.: {Multistage stochastic programs with block-separable recourse}.
\newblock Mathematical Programming Study \textbf{28}, 48--62 (1986).
\newblock \doi{10.1007/BFB0121125}

\bibitem{Mazzi2020}
Mazzi, N., Grothey, A., McKinnon, K., Sugishita, N.: {Benders decomposition with adaptive oracles for large scale optimization}.
\newblock Mathematical Programming Computation \textbf{13}, 683--703 (2020).
\newblock \doi{10.1007/s12532-020-00197-0}

\bibitem{Papavasiliou2018ApplicationUncertainty}
Papavasiliou, A., Mou, Y., Cambier, L., Scieur, D.: {Application of Stochastic Dual Dynamic Programming to the Real-Time Dispatch of Storage under Renewable Supply Uncertainty}.
\newblock IEEE Transactions on Sustainable Energy \textbf{9}(2), 547--558 (2018).
\newblock \doi{10.1109/TSTE.2017.2748463}

\bibitem{pereira1991}
Pereira, M.V., Pinto, L.M.: {Multi-stage stochastic optimization applied to energy planning}.
\newblock Mathematical Programming \textbf{52}(1-3), 359--375 (1991).
\newblock \doi{10.1007/BF01582895/METRICS}

\bibitem{Shapiro2011AnalysisMethod}
Shapiro, A.: {Analysis of stochastic dual dynamic programming method}.
\newblock European Journal of Operational Research \textbf{209}(1), 63--72 (2011).
\newblock \doi{10.1016/J.EJOR.2010.08.007}

\bibitem{VanSlyke1969}
Van~Slyke, R.M., Wets, R.: {Optimal control and stochastic programming}.
\newblock SIAM journal on applied mathematics \textbf{17}(4), 638--663 (1969).
\newblock \doi{10.1137/0117061}

\bibitem{Zakeri}
Zakeri, G., Philpott, A., Ryan, D.: {Inexact cuts in Benders decomposition}.
\newblock SIAM Journal on Optimization \textbf{10}, 643--657 (2000).
\newblock \doi{10.1137/S1052623497318700}

\bibitem{Zhang2025IntegratedDecomposition}
Zhang, H., Grossmann, I.E., McKinnon, K., Knudsen, B.R., Nava, R.G., Tomasgard, A.: {Integrated investment, retrofit and abandonment energy system planning with multi-timescale uncertainty using stabilised adaptive Benders decomposition}.
\newblock European Journal of Operational Research \textbf{325}(2), 261--280 (2025).
\newblock \doi{10.1016/J.EJOR.2025.04.005}

\bibitem{Zhang2024DecompositionProgramming}
Zhang, H., Grossmann, I.E., Tomasgard, A.: {Decomposition methods for multi-horizon stochastic programming}.
\newblock Computational Management Science 2024 21:1 \textbf{21}(1), 1--24 (2024).
\newblock \doi{10.1007/S10287-024-00509-Y}

\bibitem{Zhang2024AUncertainty}
Zhang, H., Mazzi, N., McKinnon, K., Nava, R.G., Tomasgard, A.: {A stabilised Benders decomposition with adaptive oracles for large-scale stochastic programming with short-term and long-term uncertainty}.
\newblock Computers {\&} Operations Research \textbf{167}, 106665 (2024).
\newblock \doi{10.1016/J.COR.2024.106665}

\bibitem{Zou2019StochasticProgramming}
Zou, J., Ahmed, S., Sun, X.A.: {Stochastic dual dynamic integer programming}.
\newblock Mathematical Programming \textbf{175}(1-2), 461--502 (2019).
\newblock \doi{10.1007/S10107-018-1249-5/TABLES/5}

\end{thebibliography}

\end{document}